\documentclass{article}
\usepackage[T1]{fontenc}
\usepackage[english]{babel}
\usepackage{amsfonts,amsmath,amssymb,amsthm}
\usepackage{graphicx,mathtools}
\usepackage{tikz}
\usetikzlibrary{backgrounds,calc}

\newcommand{\isga}{isometric group action}

\newcommand{\pig}{principal isotropy group}

\newcommand{\CC}{\mathbb{C}}
\newcommand{\DD}{\mathbb{D}}

\newcommand{\NN}{\mathbb{N}}

\newcommand{\RR}{\mathbb{R}}
\newcommand{\ZZ}{\mathbb{Z}}

\newcommand\PP{\mathbb{P}}

\newcommand{\D}{\operatorname{D}}
\newcommand{\G}{G}
\renewcommand{\H}{H}
\newcommand{\K}{K}
\newcommand{\N}{N}
\newcommand{\Pii}{\operatorname{\Pi}}
\newcommand{\Pin}{\operatorname{Pin}}

\renewcommand{\S}{\operatorname{S}}

\newcommand{\SO}{\operatorname{SO}}
\newcommand{\Sp}{\operatorname{Sp}}
\newcommand{\SU}{\operatorname{SU}}
\newcommand{\T}{\operatorname{T}}
\newcommand{\U}{\operatorname{U}}

\newcommand{\p}{\operatorname{p}}

\newcommand{\M}{\mathcal{M}}
\newcommand{\TT}{\mathcal{T}}
\newcommand{\B}{\mathcal{B}}
\newcommand{\W}{\mathcal{W}}

\newtheorem{thm}[equation]{Theorem}
\newtheorem{cor}[equation]{Corollary}
\newtheorem{lem}[equation]{Lemma}
\newtheorem{prop}[equation]{Proposition}
\newtheorem*{thm*}{Theorem}
\newtheorem*{cor*}{Corollary}
\newtheorem*{lem*}{Lemma}
\newtheorem*{prop*}{Proposition}
\newtheorem{THM}{Theorem}[]

\theoremstyle{definition}
\newtheorem{defn}[equation]{Definition}
\newtheorem*{defn*}{Definition}
\theoremstyle{remark}
\newtheorem{rem}[equation]{Remark}
\newtheorem{ex}[equation]{Example}

\newtheorem*{rem*}{Remark}
\newtheorem*{ex*}{Example}
\newtheorem*{assumption*}{Assumption}
\newtheorem*{remnot}{Notation}

\newtheorem*{pf}{Proof}

\newcommand{\cref}[1]{Corollary~\ref{#1}}
\newcommand{\dref}[1]{Definition~\ref{#1}}
\newcommand{\eref}[1]{Example~\ref{#1}}
\newcommand{\lref}[1]{Lemma~\ref{#1}}
\newcommand{\pref}[1]{Proposition~\ref{#1}}
\newcommand{\rref}[1]{Remark~\ref{#1}}

\newcommand{\tref}[1]{Theorem~\ref{#1}}
\newcommand{\fref}[1]{Figure~\ref{#1}}
\numberwithin{equation}{section}

\sffamily

\newcommand{\sect}[1]{\section{\mdseries{#1}}}

\begin{document}

\title{\textsc{
\textbf{
Low Dimensional Polar Actions
}}}
\date{}
\author{Francisco J. Gozzi\footnote{The author was supported by CNPq-BRAZIL.}}
\maketitle


The concept of symmetry has played an important role in geometry since its inception.
Within the field of Riemannian geometry we find the general theory of proper isometric actions 
by Lie groups, and as a special case of interest that of polar actions.

A complete Riemannian manifold $M$ together with a proper isometric action of a Lie group $\G$ 
is said to be \textit{polar} if it admits a \textit{section}, i.e., an immersed complete submanifold 
$\Sigma$ of $M$ intersecting every $\G$-orbit and doing so orthogonally. 
The existence of $\Sigma$ is equivalent to the integrability of the horizontal distribution along the regular part;
see \cite{Boaulem}.  
The section is preserved by the so called polar group $\Pii$ which is a discrete group acting properly and discontinuously \mbox{on $\Sigma$}
with quotient $\Sigma/\Pii$ isometric to $M/\G$. Hence, the orbit space is a good orbifold. 
The existence of a section enables the reduction of a geometric problem on $M$ invariant under $\G$ to a possibly simpler 
one on $\Sigma$ invariant \mbox{under $\Pi$}, as e.g. 
illustrated by the generalization of the Chevalley Restriction Theorem in \cite{PalaisTerng}
which states that smooth $\G$-invariant functions on $M$ are in one to one correspondence with $\Pii$-invariant functions on $\Sigma$,
$C^{\infty}(M)^{\G}\cong C^{\infty}(\Sigma)^{\Pii}$.

The knowledge of the orbit space together with certain isotropy data 
is enough to reconstruct the equivariant diffeomorphism type of the action, as shown in \cite{GZ}. 
This feature provides flexibility in the construction of new examples in terms of abstract polar data
and facilitates their equivariant classification. 
It also proves useful in adapting general tools, as e.g. equivariant surgeries, to the polar context.

Among classical examples are the adjoint action of a compact Lie group on itself, equipped with 
a bi-invariant metric, or more generally 
the left action \mbox{of $\K$} on $\G/\K$, associated to a symmetric pair $(\G,\K)$. 
The isotropy representations of 
$\K$ on $T_{[e\K]} \G/\K$, so called \textit{s-representations},
are of particular importance since they agree with the class of linear polar actions up to orbit equivalence; 
\mbox{see \cite{Dadok}}. In fact, such $s$-representations completely describe the local structure of a general polar action.

Trivial examples are given by the product of a homogeneous space with a manifold.
Another simple yet rich source of examples is that of cohomogeneity one actions, 
where a section is given by a geodesic orthogonal to a regular orbit. 
In this case, polar data is given by a group diagram \mbox{$ \H \subseteq \K_-,\K_+ \subseteq \G $}, 
and the reconstruction of the manifold is achieved as a gluing of two disc bundles 
$\G \times_{\K_{\pm}}\DD^{l_{\pm}}$ along their common boundary $\G/\H$.
In this framework, compact simply-connected cohomogeneity one manifolds were classified up to \mbox{dimension $7$} in \cite{Hoelscher}. 
In particular, in dimension $5$ the author shows that 
$M$ is diffeomorphic to either $\mathbb{S}^5$, $\SU(3)/\SO(3)$, $\mathbb{S}^3\times \mathbb{S}^2$ or the 
non-trivial $\mathbb{S}^3$-bundle over $\mathbb{S}^2$, denoted by $\mathbb{S}^3\tilde{\times} \mathbb{S}^2$. In addition, he classifies the actions up to equivariant diffeomorphism.

\smallskip
In this paper we classify all closed simply-connected polar manifolds of dimension $5$ or less 
by addressing the case of cohomogeneity at least two.
We assume that the group acting is connected, since a polar action by a \mbox{Lie group $\G$}	 
naturally restricts to its identity component $\G_0$ in a polar fashion and with the same cohomogeneity.
In order to facilitate the description of the actions we shall assume they are almost-effective, i.e., at most a finite
subgroup acts trivially.
Polar manifolds of dimension $3$ or $4$ are not hard to classify (see Theorems \ref{classification dim3} and \ref{classification dim4})
where in dimension $4$ the main case of $\T^2$ actions follows from \cite{OR,GZ}.
Our main contribution is thus the classification of polar \mbox{$5$-manifolds}. 
Among these, the effective polar $\T^2$-actions are the most important and difficult case.
We describe both the equivariant type of the actions as well as the diffeomorphism type of the underlying manifolds.

Many of the actions discussed later on can be described as follows.

\smallskip
\noindent\textbf{Main Example.}
Consider $\mathbb{S}^3\times \mathbb{S}^3 \subset \CC^4$ where the circle acts linearly by multiplication as
$\theta \star (z_1,z_2,z_3,z_4)= (z_1\theta,z_2\theta^{-1},z_3\theta^k ,z_4)$, $\theta \in \CC, |\theta|=1$, $k\in \ZZ$. 
We denote this circle by $\S^1_{(1,-1,k,0)}$ and observe that it acts freely on $\mathbb{S}^3\times\mathbb{S}^3$.
Passing to the quotient, we retrieve one of the 
two $\mathbb{S}^3$ bundles over $\mathbb{S}^2$ with map to the base induced from the projection onto the first factor of $\mathbb{S}^3\times \mathbb{S}^3$.
In fact, the outcome manifold depends on the parity of $k$, as we may trivialize such bundle over two hemispheres $\DD^2_{\pm}\subset \mathbb{S}^2$ 
and observe that the gluing along their common boundary $\mathbb{S}^1$ is determined by $k \in \pi_1(\SO(3))\cong \ZZ_2$. We have:
$$\mathbb{S}^3\times_k \mathbb{S}^2:= (\mathbb{S}^3\times\mathbb{S}^3)/\S^1_{(1,-1,k,0)} \cong 
 \begin{cases}
  \mathbb{S}^3 \times 		\mathbb{S}^2 	&	k ~ \text{even}\\ 
  \mathbb{S}^3 ~\!\tilde{\times} ~\mathbb{S}^2  	&	k ~ \text{odd}.	
  \end{cases} $$
We may construct natural actions on the quotients $\mathbb{S}^3\times_k \mathbb{S}^2$, induced from actions on $\mathbb{S}^3\times \mathbb{S}^3$ that 
either extend or commute with the previous right action by $\S^1_{(1,-1,k,0)}$: 
\begin{enumerate}

\item A circle action by $\S^1_{(0,0,0,1)}$ acting as $\theta\star (z_1,z_2,z_3,z_4)=(z_1,z_2,z_3,z_4\theta)$;

\item The standard action of $\SU(2)$ on the first factor of $\mathbb{S}^3\times \mathbb{S}^3$;

\item An $\SU(2)\times \S^1$ action 
by $(A,\theta)\star (z_1,z_2,z_3,z_4)= (A \cdot(z_1,z_2),z_3\theta,z_4)$;

\item A linear $\T^3$ action by $(\theta_1,\theta_2,\theta_3,1)\star (z_1,z_2,z_3,z_4)= (z_1\theta_1,z_2\theta_2,z_3\theta_3,z_4)$,
which induces an effective action by $\T^2\cong \T^3/\S^1_{(1,-1,k,0)}$ on $\mathbb{S}^3\times_k \mathbb{S}^2$;
\item The maximal torus $\T^4\subset \SO(4)\times \SO(4)$ which induces an effective $\T^3$ action on $\mathbb{S}^3\times_k \mathbb{S}^2$. 
\end{enumerate}
As we will see, the induced actions 
are all polar with respect to the submersion metric on $\mathbb{S}^3\times_k\mathbb{S}^2$, 
inherited from the standard metric on $\mathbb{S}^3(1)\times \mathbb{S}^3(1)$. In particular, the metrics have non-negative curvature.

\smallskip
We now proceed with the description of all polar $\T^2$ actions on $5$-manifolds,
where the above $\T^2$ actions on $\mathbb{S}^3\times_k \mathbb{S}^2$
play a key role as basic building blocks of our construction.
\begin{THM} \label{THM T2 in M5}
 A closed simply-connected $5$-manifold equipped with an effective polar action by $\T^2$ is obtained,
by equivariant surgery operations, from linear polar actions on spheres and products of spheres,
or one of the above $\T^2$-actions on $\mathbb{S}^3\times_k \mathbb{S}^3$.
Moreover, 
the equivariant surgeries correspond either to connected sums at fixed points, or surgeries along regular orbits. 
\end{THM}
The question as to which surgeries are allowed is more subtle, see \fref{fig surgeries} 
for some examples and Section \ref{section Polar $T^2$ actions on $M^5$} for a detailed discussion.

For the diffeomorphism type of the manifolds themselves, we use the Barden-Smale classification of $5$-manifolds.
It turns out that for polar $\S^1$, $\T^2$, and $\T^3$ actions the possible diffeomorphism types agree.
\begin{THM} \label{THM diffeo tori in M5}
A closed simply-connected $5$-manifold that admits a polar 
action by an abelian Lie group is diffeomorphic to either $\mathbb{S}^5$, $\mathbb{S}^3\times \mathbb{S}^2$, $\mathbb{S}^3\tilde{\times} \mathbb{S}^2$, 
or their connected sums.
\end{THM}

The case of $\T^3$ actions in dimension $5$ is treated in \cite{Oh5}. 
More generally, all cohomogeneity two tori actions are automatically polar, as observed in \cite{GZ}.
The case of polar $\S^1$ actions is fairly simple, and can also be treated in general for any dimension; 
see \pref{cohomogeneity n-1}. 

For the classification of the non-abelian case we have:
\begin{THM} \label{THM classification dim5} 
Let $\G$ be a compact connected non-abelian Lie group, with a non-trivial and almost-effective polar action on a compact simply-connected \mbox{$5$-manifold} $M$
of cohomogeneity at least two. 
Then $M$ is equivariantly diffeomorphic to one of the following:
\begin{enumerate}
\item[a)] The linear polar actions on $\mathbb{S}^5$; 
\item[b)] The $\SU(2)$ and $\SU(2)\times \S^1$ actions on $\mathbb{S}^3 \times_k\mathbb{S}^2$ as described in the Main Example;
\item[c)] The factor-wise linear actions by $\SO(3)\times \S^1$ and $\SO(3)\times \{1\}$ on $\mathbb{S}^3 \times \mathbb{S}^2$, 
and their connected fixed-point sums $\#_{n} \mathbb{S}^3\times \mathbb{S}^2$;
\item[d)] The unique $\SO(3)$-actions on the Wu manifold $\W=\SU(3)/\SO(3)$ and on 
the Brieskorn variety $\B$ of type $(2,3,3,3)$, 
or their respective fixed-point connected sums, $\#_{k}\B$ and $\#_{l}\W$.
\end{enumerate}
\end{THM}
In \cite{Simas} general actions by $\SU(2)$ and $\SO(3)$ were classified. 
It follows from our work that all but 2 of the actions admitting singular orbits are in fact polar.

We will see that further fixed point connected sums $\#_{k} \B \#_{l}\W$ are equivariantly diffeomorphic to $\#_{k}\B$ or $\#_{l}\W$
themselves, and that $\#_{k}\B$ has a unique polar action by $\SO(3)$ whereas $\#_{l}\W$ has finitely many.

As an immediate consequence of the previous theorems we have the following,
\begin{cor*}
 A closed simply-connected polar $5$-manifold is diffeomorphic to either
 $\mathbb{S}^5$, $\#_{\alpha}(\mathbb{S}^3\tilde{\times} \mathbb{S}^2)\#_{ n_0}(\mathbb{S}^3\times \mathbb{S}^2)$, $\#_{n_1}\B$ or $\#_{n_2}\W$, for $n_i \in \NN_0$ and $\alpha=0,1$.
 \end{cor*}

As an application, we obtain an equivariant classification of all polar actions on \mbox{$5$-manifolds} which admit invariant metrics with non-negative curvature. 
\begin{THM} \label{THM non-negative curv}
 Let $M$ be a closed simply-connected $5$-manifold with a polar action of cohomogeneity at least two 
 whose sectional curvature is non-negative. Then $M^5$ is equivariantly diffeomorphic to 
a product of linear actions on spheres, 
the left action of $\SO(3)$ on $\SU(3)/\SO(3)$, or one of the actions on $\mathbb{S}^3\times_k \mathbb{S}^2$ constructed in the Main Example above.
\end{THM}
In contrast, cohomogeneity one actions on $5$-manifolds all admit non-negative curvature, see \cite{Hoelscher}.
The added rigidity when the cohomogeneity is $\geq 2$ is analogous to the case of positive curvature, where it was shown that 
a polar action is equivariantly diffeomorphic to a linear action on a compact rank one symmetric space; see \cite{FGT}.

\smallskip
The paper is organized as follows.
In Section \ref{section Structure of polar manifolds} we summarize the basic definitions and results on polar actions
and describe the equivariant surgeries on such actions. 
The second section is devoted to the construction of various examples that will appear throughout the paper.
In Section \ref{section Classification in low dimensions}, we discuss all actions other than the ones by $\T^2$ in dimension $5$, 
which are treated in \mbox{Section \ref{section Polar $T^2$ actions on $M^5$}}. The proof of \tref{THM diffeo tori in M5} is the content 
of \mbox{Section \ref{sect diffeo tori classif}}, while the last section is reserved for applications in non-negative curvature.

\medskip
This work was completed as part of the author's Ph.D. thesis at IMPA. 
He would like to thank the institution for its hospitality and support, 
and express his deepest thanks to his doctoral supervisors 
Luis A. Florit and Wolfgang Ziller 
for their friendly advice and patient guidance.

\sect{Preliminaries on Polar Manifolds \label{section Structure of polar manifolds}}
 
We begin by recalling general facts about isometric group actions (see e.g. \cite{Bredon}) and fixing some notations. 
Next, we discuss polar manifolds as described in \cite{GZ}, together with the equivariant ``cut and paste'' 
operations needed for the construction of large families of relevant examples.

\medskip
Let $\G$ denote a compact Lie group acting by isometries on a Riemannian manifold $M$. 
The orbit through a point $p\in M$ 
is an embedded homogeneous submanifold $ \G\cdot p \cong \G /\G_p$,
where $\G_p=\{g\in \G \mid g.p=p \}$ is the \textit{isotropy subgroup} at $p$. 
This group acts on the normal sphere to the orbit in $T_p(\G\cdot p)^{\perp}$, and by the Slice Theorem we have: 
\begin{equation} \label{slice theorem}
\TT_{\epsilon}(\G\cdot p) \cong (\G \times S_p)/ {\G_p}=: \G \times_{\G_p} S_p,
\end{equation}
where $S_p$ 
is a \textit{slice} at $p$. 
The neighborhood of a point in the orbit space, \mbox{$[p] \in M/\G$}, is homeomorphic to the quotient $S_p/\G_p$.
The conjugacy class $(\G_p)$ of an isotropy subgroup is called an \textit{isotropy type}.
Orbits corresponding to a minimal isotropy type are called \textit{regular}, and we reserve $\H$ to denote a choice of
one of this regular or \textit{principal} isotropy subgroups.
Orbits with bigger isotropy type are called \textit{singular},
unless they have the same dimension as the principal ones in which case they are called \textit{exceptional}.

The orbit space is not a manifold in general but a disjoint collection of strata,
which are smooth manifolds corresponding to the projection of each $M_{(\K)}:= \{x\in M ~|~ (\G_x)=(\K)\}$. 
The partial order between isotropy types induces an inverse relation between the strata,
where a smaller isotropy type $(\K_1)\leq (\K_2)$ corresponds to a bigger stratum containing the smaller one in its closure, 
$M_{(\K_2)}\subset \overline{M_{(\K_1)}}$. 
%
%
%
%

\smallskip
In the case of polar actions we obtain an optimal reduction in dimension for the understanding of the orbit space $M/\G$ 
and its stratification by isotropy types. 
In fact, the action of $\G$ on $M$ induces an action of a discrete group 
$\Pii=N(\Sigma)/ Z(\Sigma)$, called the \textit{polar group}, where 
$N(\Sigma)=\{g\in \G | g\cdot\Sigma=\Sigma\}$ 
and the subgroup $Z(\Sigma)$ is the point-wise stabilizer of the section. 
The latter agrees with the \pig~of a regular point in the section, i.e., $Z(\Sigma)= \H$, while 
$\Pii$ acts properly and discontinuously on 
$\Sigma$ with isometric quotient $\Sigma/\Pii \cong M/\G$.
In this way, the section provides an orbifold cover of the orbit space.

Inside $\Pii$ we identify \textit{reflections}, i.e., elements which have a hypersurface or \textit{mirror} 
$\M^{k-1} \subset \Sigma^k$ inside their fixed-point set.
We consider the complement \mbox{in $\Sigma$} of all possible mirrors $\{\M_i\}$ and define a \textit{chamber} $C$
to be the closure of a connected component.  
By construction, the boundary of $C \subseteq \Sigma$ has a stratification by totally geodesic faces, given by the
mirrors and their intersections. 	
More precisely, we define \textit{faces}
$F_i$ of $C$ to be the closure of each connected component of the open intersections with a mirror, 
$$F_i=\overline{ int(\M_i \cap C)} \subset \M_i.$$
Deeper faces $F_{i_1\cdots i_k}$ 
are defined as the connected intersections of the previous ones,
$F_{i_1\cdots i_k} = F_{i_1}\cap \cdots \cap F_{i_k}$ 
The length of the sub-index 
is referred to as the \textit{depth} at which the face lies. To each face we can assign a unique generic isotropy subgroup.


Recall that by definition every orbit intersects the section, and thus a given chamber as well.
Moreover, an orbit may possibly intersect $C$ more than once, implying the existence of a non trivial 
stabilizer subgroup preserving the chamber, \mbox{$\Pii_C\!:=\{g\in \Pii| g\cdot C =C\}\subset \Pii$}.

From now on, we shall restrict ourselves to the following class of actions, which is a natural setting in our context 
since a polar action by a connected group on a simply-connected manifold is always of this type (see \rref{coxeter is natural}). 
\begin{defn} \label{coxeter definition}
A polar action without exceptional orbits and with trivial stabilizer group $\Pii_C$ is called \textit{Coxeter polar}.
\end{defn}
One of the main properties about Coxeter polar actions is that the inclusion $C \to \Sigma \xrightarrow{\sigma} M$ 
induces an isometry between the orbit spaces and the chamber,
\begin{equation} \label{orbit space}
C \cong \Sigma /\Pii \cong M/\G.
\end{equation}
In the case of a general \isga~we may not have a cross section to the projection map $M \to M/\G$, 
while for a Coxeter polar action this is given by the chamber $C$. 
Furthermore, the faces have constant isotropy group along the interior, and 
\eqref{orbit space} identifies the open strata in $C$ with the $\Pii$-isotropy type strata in $\Sigma/\Pii$, 
and the corresponding $\G$-isotropy type strata in $M/\G$. 
We may encode this information in a marking of the open strata $F_I$ of $C$ 
forming a partially ordered graph of $\G$-isotropy subgroups $\K_I$ associated to each one.
The subgraph below a given vertex $K_I$ is called the \textit{history} of $K_I$ and corresponds to the labeling of strata at a point $p\in F_I$.
We say that the marking of $C$ is \textit{compatible} if the history at a point $p \in F_I$
comes from a Coxeter polar representation, 
and the orbifold group at $p\in C$ coincides with the polar group of such representation. 


The isotropy subgroups of $\Pii$ are Coxeter groups, in fact given $p \in F_{i_1\cdots i_k}$, 
we have that $\Pii_p$ is generated by the reflections $r_i$ defined by the corresponding
mirrors through $p$, and satisfying relations $(r_i.r_j)^{m_{ij}}=1$ depending on the angle between them. 
The whole polar group $\Pii$ is a quotient of a Coxeter group, as it is generated by the reflections along mirrors in $\Sigma$, but 
there may be more relations than the local ones just described.
We shall regard the orbit space as a \textit{Riemannian Coxeter orbifold}, 
i.e., a Riemannian orbifold modeled on finite reflection groups. 
It is a good orbifold, with only manifold points along the interior. 
\begin{defn} \label{polar data}
A connected Riemannian Coxeter orbifold $C$ 
together with a compatible group graph marking the strata defines \textit{Coxeter polar data}, denoted by $(C,\G(C))$.
\end{defn}


We point out that the slice representations do not need to be specified as part of the Coxeter polar data. 
Indeed, they are determined by the group marking due to the following lemma that will be proved in general in a forthcoming paper. 
Observe that this is a weaker statement than that of item (b) of \mbox{Lemma 2.4} in \cite{GZ}, 
whose proof seems to have a gap for reducible representations. 
An ad-hoc argument for representations of dimension five or less follows easily from the work of \cite{Bergmann,EH_classif_polarrep}
and would also be enough for the scope of this paper.

\begin{lem}
 \label{coxeter rigidity}  
A Coxeter polar representation by a 
Lie group $L$ is determined up to linear equivalence by 
its history and dimension. Furthermore, $L$ is generated by the face isotropy groups, as well as by $L_0$ and $\H$. 
 \end{lem}

In particular, the group graph of a Coxeter polar action must satisfy the following algebraic conditions: 
\begin{equation} \label{coxeter rigidity2}
\begin{array}{c}
 \K_{i_0 \cdots i_l}	= \text{gen}_{0\leq k\leq l} \{\K_{i_k}\}	\\
									\\
 \H 		=\K_{i} \cap \K_j,~\forall i\neq j \text{ with } F_{ij}\neq \emptyset.
\end{array}
\end{equation}
The first equation rephrases part of the previous Lemma. 
The second one is valid for adjacent faces, i.e., those with non-empty intersection. 
In fact, by definition $\K_i \supset \H$, 
and $\K_i \cap \K_j \subset \H$ since it acts on the slice at $p\in F_{ij}$ with fixed point set 
including the tangent space to each face, and therefore spanning all of $T_p\Sigma$. 

The face isotropy group $\K_i$ acts
transitively on the sphere orthogonal to its strata, 
\begin{equation} \label{next to minimal}
\K_i/\H \cong \mathbb{S}^{l_i}, ~ l_i > 0. 
\end{equation}
Because of \eqref{coxeter rigidity2}
the group graph is completely determined by the isotropy groups $\K_i$ marking the maximal faces. 
Furthermore, the reflections $r_i \in \Pii$ along these maximal faces generate the polar group. 
These elements can also be thought of as the unique involutions in the normalizer $\N^{\K}(\H)/\H \subset \N(\H)/\H$. 

As an application of the fact that $\Pii$ acts simply-transitively on the set of chambers, 
thus \textit{tessellating} the section, we can try to orient $\Sigma$ by prescribing 
an orientation on $C$ and declaring each reflection along a maximal face to be orientation reversing.
The reader may see \ref{ex T3 case main example} for an illustrative example having non-orientable section.
This leads to the following criterion for orientability.  
\begin{lem} \label{lem section orientability}
The section of a Coxeter polar manifold is orientable if and only if the chamber is, and if there are no relations
of odd order among the reflections corresponding to maximal faces.
\end{lem}

%
%
%

Coxeter polar data 
is an invariant of the action which is enough to construct and identify Coxeter polar manifolds: 
\begin{thm}[\cite{GZ}] \label{GZ coxeter reconstruction} 
Coxeter polar data $(C,\G(C))$ determines a Coxeter polar manifold $M(C,\G(C))$ up to equivariant diffeomorphism. 
\end{thm}
We point out though that 2 different Coxeter polar data $(C,\G(C))$ and $(C,\G'(C))$
can be equivariantly diffeomorphic since the diffeomorphism may not respect the metric. 
This can also change the polar group and section, in fact, a closed section can become non compact. 
It is an interesting question whether via such a change of metric the section for any Coxeter polar action can be made compact, 
as is the case for cohomogeneity one manifolds.
See \cite{GWZ} for a discussion in the case of cohomogeneity one actions.

As an immediate application of \tref{GZ coxeter reconstruction} observe that each strata in $M$ gets uniquely reconstructed as a trivial bundle 
$M_{(\K)}= M_{(\K)}^*\times \G\!/\!\K $,
since it corresponds to data given by the manifold $M_{(\K)}^*$ constantly marked by $\K$.

\begin{rem} \label{coxeter is natural}
A general polar manifold can be described as the quotient of a Coxeter polar manifold by a group acting freely and discontinuously which
commutes with the action and preserves orbit types; see 
\cite{GZ}. 
In particular, it follows that a simply-connected manifold with a polar action by a connected group is Coxeter polar. 
\end{rem}
The following lemmas show how to construct new polar actions as quotients of given ones. 

\begin{lem} \label{lem commuting perp action is polar}
Let $M$ be a polar $\G$-manifold, and $L$ be a Lie group acting freely on $M$, such that the two actions commute.
Then the induced $\G$ action on $M/L$ is polar if $\G$ and $L$ have orthogonal orbits, 
i.e., $T_p(\G\!\cdot p) \perp T_p (L\!\cdot p)$.
\end{lem}
\begin{pf} The assumption implies that $T_p (L\!\cdot p)\subset T_p \Sigma$ and hence  $\Sigma$ is $L$-invariant. Furthermore, 
it projects to an immersed complete submanifold in $M/L$. 
The polarity of the action is clear once we observe that $\Sigma/L$ has the right dimension and 
integrates the horizontal distribution along the regular part. 
\qed
\end{pf}
%
\begin{lem} \label{lem normal subgroup quotient action is polar}
Let $M$ be a (Coxeter) polar $\G$-manifold with section $\Sigma$, and $L \triangleleft \G$ be a normal subgroup which acts freely on $M$.
Then the induced $\G/L$ action on the quotient $M/L$ is (Coxeter) polar with the submersion metric on $M/L$. A section is given
as a proper and discontinuous quotient of $\Sigma$ by the action of $\Pii \cap L$.
Moreover, in the Coxeter polar case, chambers are isometrically identified and the marking by $\G\!/\!L$-isotropy groups is
obtained as the projection along $\G\to \G\!/\!L$ of the previous marking of $\G$-isotropy groups. 
This gives an isomorphism of isotropy groups.
In addition, the new polar group is obtained as the quotient $\Pii / \Pii\cap L$.
\end{lem}
\begin{pf}
The section $\sigma:\Sigma\to M$ is horizontal with respect to the projection map	 $\p:M\to M/L$, 
therefore we have by composition a section in the quotient given by $\p\circ\,\sigma\!:\!\Sigma\to M/L$. 
The orbit spaces are isometric to each other, $M/\G=(M/L)/(\G/L)$.
It then follows, that in the Coxeter polar case $\p$ identifies chambers isometrically. 
The isotropy group at a point $\p(x) \in M/L$ corresponds to the projection of 
$G_x$ under the map $\Phi:\G\to \G\!/\!L$. In fact, the projection $\Phi$ 
induces an isomorphism between isotropy groups since $L\cap G_x= \{e\}$.
This shows that we have retrieved a Coxeter polar action associated to data $(C,\G\!/\!L(C))$. 

Notice in general that we have a unique minimal section up to covers, which admits possible self-intersections but does not have
any self-tangencies since it is totally geodesic. 
In our case, $\Sigma$ covers the new section in $M/L$ up to the discrete action of $\N(\Sigma) \cap L$. 
As before, we may identify $\N(\Sigma)\cap L=:\Pii \cap L$ since $L\cap \H=\{e\}$.
Finally, in the Coxeter polar case, the polar group is generated by the same reflections (as elements in $\G$) 
thanks to the coincidence of the chambers, but is possibly smaller and equals $\Pii/\Pii\cap L$. 
\qed
\end{pf}
%

\begin{rem} \label{imprimitive}
Given a (Coxeter) polar $L$-action on $M'$, and $L\subset \G$, we can construct a 
(Coxeter) polar $\G$-manifold $M$ as the quotient of $M' \times \G$ 
by the diagonal action of $L$, 
 \begin{equation} 
  M \cong M' \times_{L} \G,
 \end{equation}
where the action of $\G$ is given by right multiplication on the second factor. 
Furthermore, the inclusion of a section $\Sigma'$ of $M'$ as $\Sigma\times \{e\} \to M' \times_{L} \G$ 
gives a section for the $\G$-action on $M$.
In particular, the action on $M$ is \textit{impritive}, i.e., it admits an equivariant map $M\to \G/L$ for $L\subsetneq \G$.

Conversely, consider a Coxeter polar $\G$-manifold $M$ such that all the isotropy subgroups in the data $\G(C)$ 
are contained in a closed subgroup $L\subseteq \G$. 
Then the same Coxeter polar data yields, by restriction, compatible data $(C,L(C))$ for an $L$-manifold $M'$. 
By uniqueness of the reconstruction process we have that $M$ agrees with $M' \times_{L} \G$ as above,
up to an equivariant diffeomorphism.	
\end{rem}

\medskip

Finally, there are equivariant ``\textit{cut and paste}'' operations on $\G$-manifolds that can be adapted to work on polar manifolds.
Given two $\G$-manifolds $M_1$, $M_2$, we may glue them along orbits of the same \textit{slice type}, i.e.,
same isotropy type and isomorphic slice representations.
In this case, orbits will not only be diffeomorphic,
but will also have equivariantly isomorphic tubular neighborhoods $\TT_1:=\TT(\G\cdot {p_1})\cong \TT(\G\cdot {p_2})=:\TT_2$.
We can then remove these tubular neighborhoods and glue the remaining pieces along the equivariantly identified boundary. 
We thus obtain a new $\G$-manifold, denoted by 
\begin{equation} \label{orbit sum characterization}
M_1 \#_{p_1 \sim p_2} M_2: = ~M_1 \backslash \TT_1 \cup_{\partial\TT} M_2 \backslash \TT_2 . 
\end{equation}
Notice that the orbit space of $M_1 \#_{p_1 \sim p_2} M_2$ corresponds to a connected sum of $M_1/\G$ and $M_2/\G$ 
along $[p_1]\sim[p_2]$. This motivates the following. 
\begin{defn} \label{orbit sum}
Given points $p_i \in C_i\subset M_i$, $i=1,2$, lying on respective chambers of Coxeter polar $\G$-manifolds,
having the same isotropy subgroups and isomorphic slice representations,
we define their \textit{orbit sum}
as a (marking preserving) connected sum of the chambers $C:= C_1\#_{p_1\sim p_2}C_2$. For this, we  
cut out the points $p_i$ together with their neighborhoods in $C_i$, and then glue along the boundary
via the identification induced from an isomorphism between the slice representations. 
As usual, the choice of compatible orientations on the chambers determines the gluing.
Moreover, we can make the identification isometric if we first modify the metrics on each chamber, 
so that it is a product near the boundary of each ball $B_{\epsilon}(p_i)\subset C_i$. 

By \tref{GZ coxeter reconstruction} 
we have a unique Coxeter polar manifold associated to this new compatible data $(C,\G(C))$. It will be referred to
as the \textit{orbit sum} of $M_1$ and $M_2$ along the prescribed orbits, and denoted as before by $M_1 \#_{p_1 \sim p_2} M_2$. 
\end{defn}

From now on we focus on the case where the cohomogeneity is at least two. 
Taking advantage of the fact that the reconstructed manifold is not only polar but Coxeter polar, we can determine the polar group out of
the polar data $(C,\G(C))$ thanks to the fact that this group is generated by the reflections along the maximal (depth one) faces. 
Indeed, as only $0$-dimensional strata can be deleted by this operation, the assumption on the cohomogeneity implies that the gluing will not erase any maximal face,  
thus preserving the corresponding reflections. 
As a consequence we have that the new polar group $\Pii$ will be the subgroup generated by the respective polar groups 
$\Pii_i$ of the original manifolds $M_i$, when seen as subgroups of the effective normalizer to the common principal isotropy subgroup $\H$,
\begin{equation} \label{eq polar group gen}
 \Pii = gen\{ \Pii_1,\Pii_2\} \subset N(\H)/\H.
\end{equation}

Next, we discuss the two extreme cases of orbit sum, that is, along either regular orbits or fixed points. 
These will be of interest in our description of $\T^2$ polar actions on $5$-dimensional manifolds.

\begin{prop} \label{fixed point sum}
The manifold, chamber and section coming from a fixed point sum are 
the connected sum of the corresponding pieces, and the new polar group coincides with the original: 
$$ M_1 \#_{fix} M_2 	\cong  M_1 \# 	M_2 $$
$$\Sigma 		\cong \Sigma_1 \# \Sigma_2$$
$$ \Pii=\Pii_1=\Pii_2.$$
\end{prop}
\begin{pf}
Recalling the identity \eqref{orbit sum characterization} above, together with the fact that  
the tubular neighborhoods to the orbits are just balls around fixed points, 
it is clear that the resulting manifold will be the connected sum of the original ones. 

Notice that at corresponding fixed points $p_i$ we have $\Pii_i=(\Pii_i)_{p_i}$, and in turn $(\Pii_1)_{p_1}=(\Pii_2)_{p_2}$ 
as they have isomorphic slice representations. 
From \eqref{eq polar group gen} it follows that the polar group remains the same, $\Pii_1=\Pii_2=\Pii$. 

As the polar group acts freely and transitively on the set of chambers, 
we can reflect our given chamber $C= C_1\# C_2$ by the new polar group
and retrieve the section, i.e., 
$$ \Pii \cdot C = \Sigma.$$	
In particular, when reflecting the chamber minus a small neighborhood of a fixed point, we retrieve the section minus a disc $\DD^k$.
For the section of the orbit sum, we have 
\begin{equation*}
\begin{array}{lll}
\Sigma	&=& \Pii \cdot (C_1 \# C_2) ~~=~~ \Pii \cdot [(C_1 - B_{\epsilon}(p_1)) \cup_{(\partial \DD^k/\Pii)} (C_2 - B_{\epsilon}(p_2))] \\ \\
&=& [\Pii_1 \cdot (C_1 - B_{\epsilon}(p_1))] \cup_{\partial \DD^k} [\Pii_2\cdot(C_2 - B_{\epsilon}(p_2))] \\ \\
&=&(\Sigma_1-\DD^k) \cup_{\mathbb{S}^{k-1}} (\Sigma_2 - \DD^k) ~~=~~ \Sigma_1 \# \Sigma_2.
\end{array}
\end{equation*}
\qed
\end{pf}

\begin{prop} \label{regular orbit sum section}
The section associated to the polar action given by the regular orbit sum $M_1\#_{p_1\sim p_2} M_2$ 
is a connected manifold obtained by connected sums of copies of the corresponding sections $\Sigma_i$ of $M_i$, $i=1,2$,
together with additional \mbox{$0$-surgeries}.
\end{prop}

\begin{pf}
Analogously to the proof of the previous proposition, reflecting the chamber $C_1\# C_2$ retrieves the new section.

Recall that $\Pii$ is generated by the previous polar groups $\Pii_1$ and $\Pii_2$. An element $[g] \in \Pii- \Pii_1$ 
does not preserve the section $\Sigma_1$, and in fact it gives a traslated copy of it. 
Notice that the points $p_i$, along with their neighborhoods where the gluing takes place, are completely embedded in the regular part.

The spheres that bound the disks which were removed in the respective chambers $C_1$ and $C_2$, 
get reflected along with the new chamber as many as $|\!\Pii\!|$ times. 
As they lie over the regular part they do not interfere with possible self-intersections. 
Therefore, reflecting $C_1 - \epsilon.\DD^k$ by the polar group $\Pii$, will give $|\!\Pii/\Pii_1\!|$ 
disjoint copies of the original section $\Sigma_1$, with each copy 
having as many as $|\!\Pii_1\!|$ discs removed. We proceed to compute the new section:
\begin{equation*}
\begin{array}{lll}
\Sigma	&=&	\Pii \cdot (C_1 \# C_2) ~~=~~ \Pii \cdot [(C_1 - \DD^k) \cup_{\mathbb{S}^{k-1}} (C_2 - \DD^k]  \\ \\
	&=&	|\Pii/\Pii_1| \cdot [\Pii_1\cdot (C_1 - \DD^k)] 		\cup_{\Pii \cdot \mathbb{S}^{k-1}} 	|\Pii/\Pii_2| \cdot  [\Pii_2\cdot (C_2 - \DD^k)] \\ \\ 
	&=& 	|\Pii/\Pii_1| \cdot [\Sigma_1 - \Pii_1 \cdot \DD^k]	\cup_{\Pii \cdot \mathbb{S}^{k-1}}	|\Pii/\Pii_2| \cdot  [\Sigma_2 - \Pii_2 \cdot  \DD^k]\\ \\
	&=&	|\Pii/\Pii_1| \cdot \Sigma_1  ~\sqcup~ |\Pii/\Pii_2| \cdot  \Sigma_2 ~~/~|\Pii|~\text{identifications} .
	\end{array}
\end{equation*}
The removal of open discs by pairs 
for later gluing along the bounding spheres, corresponds to the definition of a $0$-surgery, which is the usual connected sum 
when the pieces are disjoint.
Keeping track of the gluings we see that the section we construct is connected, as each copy of $\Sigma_1$ 
is connected to any of $\Sigma_2$, and reciprocally. 
In particular, only $|\!\Pii/\Pii_1\!| + |\!\Pii/\Pii_2\!|-1 $ such surgeries are enough to connect the different copies. 
This gives a possibly non-trivial
remnant of $0$-surgeries to be performed on the connected sum of all the copies. 
We have:
\begin{equation*}
\begin{array}{lll}
\Sigma	&=&	 (\#_{n_1} \Sigma_1  \#_{n_2} \Sigma_2)~/~ m~~\text{\small{0-surgeries}} .
\end{array}
\end{equation*}
where $n_i = |\!\Pii/\Pii_i\!|$ and $m= |\!\Pii\!| - |\!\Pii/\Pii_1\!| - |\!\Pii/\Pii_2\!|+1$.
\qed
\end{pf}

\sect{Examples of polar actions \label{section Examples of polar actions} }

We now construct several examples that will appear in our classification. 

\begin{ex} \label{ex linear polar}
A linear polar representation gives rise to a polar action on the sphere. 
Because of \lref{coxeter rigidity} these actions are Coxeter polar whenever the group acting is connected. 
We have the following examples.
\begin{enumerate}
 \item The action of a torus $\T^n \subset \SO(2n) \subset \SO(2n+N)$ on $\mathbb{S}^{2n+N-1}\subset 
\CC^n \oplus \RR^N$ is polar, since it comes from a trivially extended product of cohomogeneity 
one circle actions of $\S^1_i$ on $\CC_i$, for $i=1,\ldots,n$. 
The section is $\{(\oplus_{i=1}^n \RR_i ) \oplus \RR^N\} \cap \mathbb{S}^{2n+N-1}$, with polar group $\ZZ_2^n$ generated by reflections 
through the origin in each $\RR_i$ factor, and corresponding chamber given by
$\{(\prod_{i=1}^n {\RR_{\geq 0}} ) \times \RR^N\}\cap \mathbb{S}^{2n+N-1}$.
Isotropy subgroups are coordinate sub-tori  $\T^l_{i_1,\dots,i_l}=\Pi_{j\in \{i_1,\dots,i_l\}} \S^1_j$. 

\item 
The irreducible $5$-dimensional $\SO(3)$ representation gives a linear action 
of cohomogeneity one on the sphere $\mathbb{S}^4$ with discrete \pig~$\H=\S(O(1)^3)\subset \SO(3)$ 
and singular isotropy groups $\K_{\pm}$ given by two different block embeddings of $O(2)$ in $\SO(3)$ in the $ij$ coordinates. 
The polar group is $\D_3$. 
Other examples are the transitive linear actions of $\SO(3)$ on $\mathbb{S}^2$ and $\SU(2)$ on $\mathbb{S}^3$, which are trivially polar.

\item Other linear polar actions on $\mathbb{S}^5$ come from the  
$\SO(3)\times \S^1$ product action on $\RR^3\oplus \RR^2$, the standard $4$-dimensional representation of $\SO(4)$, 
or its subgroup $\U(2)$, when trivially extended to $\RR^6$. The latter two are double suspensions of transitive actions on $\mathbb{S}^3$,  
and therefore have Coxeter polar data given by a $2$-disk with fixed points of the action along the boundary 
and interior marked by corresponding principal isotropy groups \mbox{$\SO(3)\subset \SO(4)$} and $\S^1\subset \U(2)$, respectively.
The first case has a bi-angle ($2$-gon) for its quotient, marked with fixed points at the vertices, sides labeled by
$K=\SO(3)$ and $\SO(2)\times \S^1$, and interior with principal isotropy group $\H=\SO(2)\subset \SO(3)$.
\end{enumerate}
 \end{ex}

Before we give more detail about the main examples presented in the introduction, let us set the following convention. 
\begin{remnot}
A circle subgroup of a given Lie group $\G$ lying inside a fixed maximal torus identified as $\T^k=\S^1\times \cdots \times \S^1$, 
will be specified by its (Lie algebra) \textit{slope} and denoted $\S^1_{{v}} \subset \T^k \subset \G$ 
where $Lie(\S^1_v)=<{v}> \subset \RR^k\cong Lie(\T^k)$, e.g. $\S^1_{(v_1,v_2)}\subset \S^1\times \S^1 \subset \G$.
Notice that ${v}$ can be chosen to have coprime integer entries. 
\end{remnot}

\smallskip
\begin{ex} \label{ex T3 case main example}
The linear $\G=\T^3$ action on $\mathbb{S}^3\times \mathbb{S}^3\subset \CC^2\oplus \CC^2$ given by 
$(\theta_1,\theta_2,\theta_3)\star(z_1,z_2,z_3,z_4)$ $= (z_1\theta_1,z_2\theta_2,z_3\theta_3,z_4)$ 
is clearly Coxeter polar since it is a product of polar actions. It admits $\mathbb{S}^1\times \mathbb{S}^2$ as section, with polar group
$\Pii\cong (\ZZ_2)^3$ given by all the involutions in $\T^3$
and chamber equal to $[-1,1] \times \DD^2$, i.e., a solid ``can''. The boundary of $C$ is composed of three faces with $\S^1$ isotropy 
and two closed circles with isotropy group isomorphic to $\T^2$. Its interior is marked by the trivial group.
Secondly, if we consider the quotient of $\mathbb{S}^3\times \mathbb{S}^3$ by the action of the circle $L=\S^1_{(1,-1,k)}$ we obtain
a Coxeter polar action on $\mathbb{S}^3\times_k\mathbb{S}^2$, as described in \lref{lem normal subgroup quotient action is polar}.
In particular the associated Coxeter polar data for both actions has the same (isometric) chamber. 
Furthermore, the marking of the induced action can be obtained as a
relabeling of the isotropy group marking given by the projection $\G \to \G/L$.
In order to facilitate its description we choose the subgroup $T^2\cong \S^1_{(0,1,0)} \times \S^1_{(0,0,-1)} \subset \T^3$, 
which is a complement to $\S^1_{(1,-1,k)} \subset \T^3$ 
and therefore projects isomorphically onto its image under the map $\Phi:\T^3\to \T^2\cong \T^3/\S^1_{(1,-1,k)}$.
We illustrate this in the following diagram:

\begin{figure}[h]
 \begin{center}
 \newcommand{\altura}{2}
\footnotesize
  \begin{tikzpicture}
\begin{scope}[rotate=-30,scale=0.9,yshift=-3cm]
 \draw[ultra thick] (0,\altura) node[right]{$\S^1_{\substack{(0,1,0)}}\!\times \S^1_{(0,0,1)}$} arc(0:360:1cm and 0.4cm) ;
\draw[dashed,thick] (0,0) arc(0:180:1cm and 0.4cm);
\draw[ultra thick] (0,0) node[right]{ $\S^1_{(1,0,0)}\!\times \S^1_{(0,0,1)}$} arc(0:-180:1cm and 0.4cm);	
\draw (0,0) .. controls +(-0.2,0.5) and +(-0.2,-0.5) .. (0,\altura)
      (-2,0) .. controls +(0.2,0.5) and +(0.2,-0.5) .. (-2,\altura);
\draw node at (-0.85,0.045) {$\S^1_{(1,0,0)}$} ;
\draw node[right] at(-1.8,0.35*\altura) {$\S^1_{(0,0,1)}$}; 
\draw node at(-0.9,\altura+0.045){$\S^1_{(0,1,0)}$};
\draw node at (0.5,-0.9) {$N^6= \mathbb{S}^3\times \mathbb{S}^3$.};
\end{scope}
\draw[thick,->] (2,-1.4)--node[below]{$\Phi$} (3,-1.4) ;

\begin{scope}[xshift=7cm]
\begin{scope}[rotate=-30,scale=0.9,yshift=-3cm]
 \draw[ultra thick] (0,\altura) arc(0:360:1cm and 0.4cm);
\draw[dashed,thick] (0,0) arc(0:180:1cm and 0.4cm);
\draw[ultra thick] (0,0) arc(0:-180:1cm and 0.4cm);	
\draw (0,0) .. controls +(-0.2,0.5) and +(-0.2,-0.5) .. (0,\altura)
      (-2,0) .. controls +(0.2,0.5) and +(0.2,-0.5) .. (-2,\altura);
\draw node at (-1,0) {$\S^1_{(1,k)}$} ;
\draw node[right] at(-1.4,0.37*\altura) {$\S^1_{(0,1)}$}; 
\draw node at(-1,\altura){$\S^1_{(1,0)}$};
\newcommand{\pos}{(-3,1)}
    \draw[color = black] \pos node[above]{$\T^2$};
    \draw[-stealth] \pos .. controls +(0.5,-0.2) and +(-0.5,+0.3)	.. (-2.05,0);
    \draw[-stealth] \pos .. controls +(0.5,-0.2) and +(-0.5,0) 	.. (-2.05,\altura);
    \draw node at (0,-1.2) {$M^5=\mathbb{S}^3\times_k \mathbb{S}^2$.};
    \end{scope}
\end{scope}
\end{tikzpicture}
 \end{center}    
\caption{\label{fig induced T2 action S3x_k S2} \footnotesize
 A $\T^3$ action on $\mathbb{S}^3\times \mathbb{S}^3$ and the induced $\T^2$ action on $\mathbb{S}^3\times_k \mathbb{S}^2$ obtained from the quotient by a circle subgroup $\S^1_{(1,-1,k)}$. 
}
\end{figure}
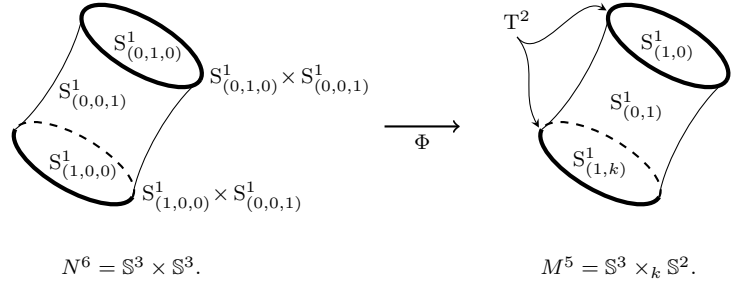

The reflections along the maximal faces for the $\T^3$ action are $(1, 1,\pm1)$, $(1,\pm 1,1)$ and $(\pm 1,1,1) \in \T^3$ 
which, by \lref{lem normal subgroup quotient action is polar}, project to reflections for the induced $\T^2$ on $\mathbb{S}^3\times_k \mathbb{S}^2$ 
given by $r_1=(1,\pm 1)$, $r_2=(\pm 1,1)$ and $r_3= (-1,(-1)^k) $, respectively. 
In addition, the new section corresponds to $\mathbb{S}^1\times \mathbb{S}^2/\tau_k$ where $\tau_k=(-1,-1,(-1)^k)$ 
is the unique involution in $L$ which belongs to the polar group of the $\T^3$ action on $\mathbb{S}^3\times \mathbb{S}^3$. 
This element reflects the $\mathbb{S}^2$ factor in $\Sigma$ along a $2$-plane when $k$ is odd and acts trivially
when $k$ is even, while it acts as the antipode along the $\mathbb{S}^1$ factor.
It is thus orientation preserving in the even case only. 
The resulting section depends on the parity of $k$:
$$ \mathbb{S}^1\times \mathbb{S}^2/\tau_k \cong \begin{cases}
                             \mathbb{S}^1\times \mathbb{S}^2 &  \text{$k$ even}\\ 
                             (\mathbb{S}^1\times \mathbb{S}^2)/\tau_k & \text{$k$ odd}. \\
                            \end{cases}$$
Notice that the group $\Pii$ is generated by the reflections $r_1,r_2$ with only even relations among them, 
but that when $k$ is odd $r_3=r_1\cdot r_2 $ gives an odd relation. 
This shows that in \lref{lem section orientability} one needs to include the reflections for all faces.

More complicated $\T^2$ actions can be constructed combining the previous ones by regular orbit sums or fixed point sums; 
see \fref{fig surgeries}. These will be discussed in more detail in Section \ref{section Polar $T^2$ actions on $M^5$}.
 \end{ex}

\begin{ex} \label{Torus cohom2}
General effective cohomogeneity two torus actions were previously studied in \cite{OR} and \cite{Oh5,Oh6}. From their work it follows that 
a $\T^n$ action on a compact simply-connected manifold $M^{n+2}$ which admits singular orbits has a topological $\DD^2$ as quotient, 
with boundary stratified by segments of isotropy $\K_i\cong \S^1$, and vertices of isotropy $\K_{i,i+1}\cong \T^2$. 
Furthermore, they show that the marking determines the equivariant diffeomorphism type of the action.
In their terminology this marking has to be \textit{legally weighted}, 
in the sense that consecutive circle isotropy subgroups give independent generators of the common vertex group $\T^2$, 
which is equivalent to \eqref{coxeter rigidity2}.
Hence, such a marking also gives compatible Coxeter polar data and hence another Coxeter polar $5$-manifold. 
From \cite{OR} it follows that these actions are equivariantly diffeomorphic. 
Thus all these manifolds turn out to be polar, as pointed out in \cite{GZ}.
%

In the case of $T^2$ actions in dimension $4$, a constructive description and classification is achieved as connected 
fixed point sums of basic actions on $\mathbb{S}^4,\pm\CC\PP^2$, $\mathbb{S}^2\times \mathbb{S}^2$ and $\CC\PP^2\# \pm \CC\PP^2$ (see \cite{OR}).
Here, the action on $\mathbb{S}^4$ is linear, while the one on $\CC\PP^2$ is obtained from the linear action of $\G=\T^3\subset \SO(6)$ on $\mathbb{S}^5$ 
under the quotient by the subgroup $L=\S^1_{(1,1,1)}$; see \lref{lem normal subgroup quotient action is polar}. 
There is also an infinite family of inequivalent actions on $\mathbb{S}^2\times \mathbb{S}^2$ and $\CC\PP^2\# -\CC\PP^2$. 
Coxeter polar data for these actions are illustrated in \fref{fig T2 in M4} below. 
\end{ex}

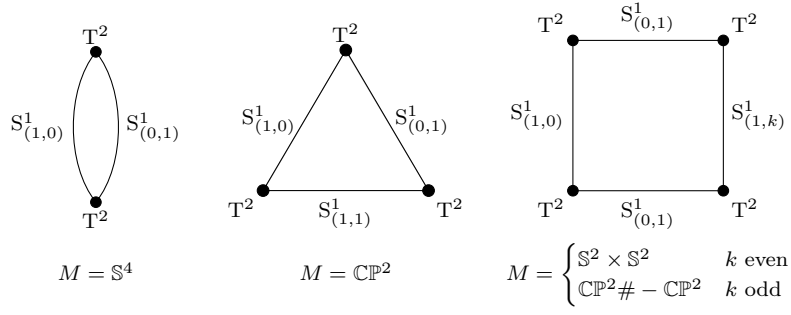
\begin{figure}[h]
\footnotesize
\begin{center}
\begin{tabular}{ccc}
\begin{tikzpicture}[scale=1]
\filldraw (0,0) circle[radius=2pt] node[above] {$\T^2$}
(0,-2) circle[radius=2pt] node[below] {$\T^2$} ;
\draw (0,0) .. controls +(0.4,-0.5) and +(0.4,0.5) .. (0,-2) node[midway,right]{$\S^1_{(0,1)}$};
\draw (0,0) .. controls +(-0.4,-0.5) and +(-0.4,0.5) .. (0,-2)node[midway,left]{$\S^1_{(1,0)}$} ;
\end{tikzpicture}  
&
\begin{tikzpicture}[scale=1.1]
\draw[fill] (0,0) circle[radius=2pt] node[below left] {$\T^2$}-- node[midway, below] {$\S^1_{(1,1)}$} (2,0) circle[radius=2pt] node[below right] {$\T^2$}
   -- node[midway, right] {$\S^1_{(0,1)}$} (1,1.7) circle[radius=2pt] node[above] {$\T^2$} -- node[midway, left] {$\S^1_{(1,0)}$} (0,0) ;
\end{tikzpicture}  
&
\begin{tikzpicture}[scale=1]

\filldraw (0,0) circle[radius=2pt] node[above left] {$\T^2$}
	  (0,-2) circle[radius=2pt] node[below left] {$\T^2$} 
	  (2,0) circle[radius=2pt] node[above right] {$\T^2$}
	  (2,-2) circle[radius=2pt] node[below right] {$\T^2$} ;
\draw (0,0) -- (2,0) node[midway, above]{$\S^1_{(0,1)}$} 
      (0,0) -- (0,-2) node[midway, left]{$\S^1_{(1,0)}$} 
      (0,-2) -- (2,-2) node[midway, below]{$\S^1_{(0,1)}$} 
      (2,0) -- (2,-2) node[midway, right]{$\S^1_{(1,k)}$} ;
	  \end{tikzpicture} 
\\
$M=\mathbb{S}^4$ & $M=\CC\PP^2$ & 
$M=\begin{cases}
\mathbb{S}^2\times \mathbb{S}^2 & k \text{ even} \\
\CC\PP^2\#-\CC\PP^2 & k \text{ odd} \end{cases}$	
\end{tabular}
\end{center}
\caption{\footnotesize \label{fig T2 in M4}
Coxeter polar data for the basic $\T^2$ actions on $M^4$. 
}
\end{figure}

In dimension $5$, there is no finite set of actions such that all $\T^3$ actions are obtained as equivariant surgeries
along orbits of $\S^1$ or $\T^2$ isotropy (\cite{Oh5}). Nevertheless, in addition to the equivariant description in terms of data, 
Oh shows that the underlying manifolds are diffeomorphic to $\mathbb{S}^5$, $\mathbb{S}^3\times \mathbb{S}^2$, $\mathbb{S}^3 \tilde{\times} \mathbb{S}^2$ or connected sums of these. 

 \begin{ex} \label{ex SU(2) of main example}
There is a polar $\SU(2)$ action on the biquotient $\mathbb{S}^3\times_k \mathbb{S}^2$.
Consider quaternion notation, $\SU(2)\cong \Sp(1)$, and let 
this group act on itself by multiplication from the left, identified as the first sphere factor of $M=\mathbb{S}^3\times \mathbb{S}^3$. 
The latter action commutes with that of $\S^1_{(1,-1,k,0)}$, 
and thus allows us to consider the product action of $\Sp(1) \times \S^1_{(1,-1,k,0)} \ni (u,\theta)$ 
on $M=\Sp(1)\times \mathbb{S}^3 $ by
$(u,\theta)\star (z_1+z_2j,z_3,z_4)=( u\cdot(z_1+z_2 j)\cdot \theta,z_3.\theta^k,z_4).$
This is easily seen to be polar with $\Sigma = \{p\}\times \mathbb{S}^2 \subset \mathbb{S}^3\times \mathbb{S}^3$, 
since it is orbit equivalent to the product of the standard action of $\SU(2)$ on the first factor with a polar circle action on the second. 
It is hence Coxeter polar by \rref{coxeter is natural}. 

We may consider the quotient under the normal subgroup $L=\S^1_{(1,-1,k,0)}$ 
to obtain an $\SU(2)$ action on $M/L=\mathbb{S}^3\times_k \mathbb{S}^2$, as in \lref{lem normal subgroup quotient action is polar}.
Coxeter polar data for the induced action is given in \fref{fig imprimitive main examples in S3xk S2} below. 
It is an imprimitive $\SU(2)$ polar action on $\mathbb{S}^3\times_k \mathbb{S}^2$.
 \end{ex}

 \begin{ex} \label{ex SU(2)xS1 of main example above}
Analogous to \eref{ex T3 case main example} is the case of the linear Coxeter polar 
$\SU(2)\times \S^1$ action on $\mathbb{S}^3\times \mathbb{S}^3$ given by
$(A,\theta)\star (z_1,z_2,z_3,z_4)= (A \cdot(z_1,z_2),z_3\theta,z_4)$ and the corresponding induced action on the quotient by 
$\S^1_{(1,-1,k,0)}$.
In particular, we retrieve an imprimitive Coxeter polar 
$\SU(2)\times \S^1$ action on the manifold \mbox{$M=\mathbb{S}^3\times_k \mathbb{S}^2$}. 
Choosing a maximal torus of $\G$ to be given by
$\U(1)\times \S^1\subset \SU(2)\times \S^1$, we have Coxeter polar data as shown in \fref{fig imprimitive main examples in S3xk S2}.
 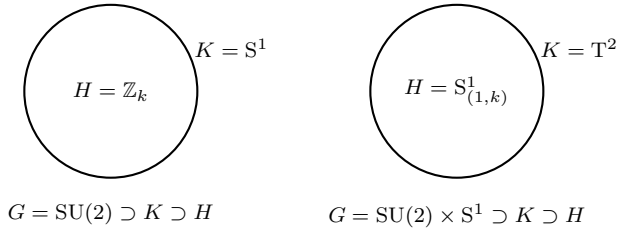
\begin{figure}[h!]
 \begin{center}
  \footnotesize  
  \begin{tikzpicture} [scale=1.15] 

\draw[thick] (0,0) node{$\H=\ZZ_k$} circle[radius=1cm] 
(0.87, 0.5) node[right]{$\K= \S^1$} 
(0,-1.2) node[below] {$\G=\SU(2)\supset \K \supset \H$} ;

\begin{scope}[xshift=4cm]
  \draw[thick] (0,0) node{$\H=\S^1_{(1,k)}$} circle[radius=1cm]
  (0.87, 0.5) node[right]{$\K= \T^2$}  
  (0,-1.2) node[below] {$\G=\SU(2)\times \S^1\supset \K \supset \H$} ;
\end{scope}
  \end{tikzpicture}
 \end{center}
 
 \caption{\label{fig imprimitive main examples in S3xk S2} \footnotesize
 Coxeter polar data for 
 actions on $M^5=\mathbb{S}^3\times_k \mathbb{S}^2$. }
\end{figure}

Notice that the diagrams in \fref{fig imprimitive main examples in S3xk S2}, 
when considered as Coxeter polar data for $\K$ actions and later 
made effective, agree with the standard polar circle action on $\mathbb{S}^3$.
\end{ex}

\begin{ex} \label{ex polar circle action non-neg}
Consider the polar action of the circle $\S^1_{(0,0,0,1)}$ on 
$\mathbb{S}^3\times \mathbb{S}^3$ which, as before, commutes with that of the circle $\S^1_{(1,-1,k,0)}$. 
The polarity of the induced $\S^1$ action on $\mathbb{S}^3\times_k \mathbb{S}^2$ follows from \lref{lem commuting perp action is polar}.
Notice that the original circle action has section $\Sigma=\mathbb{S}^3\times \mathbb{S}^2\subset \mathbb{S}^3\times \mathbb{S}^3$ with chamber $\mathbb{S}^3\times \DD^2$,
and that in this case the action of
$\S^1_{(1,-1,k,0)}$ restricts to $\Sigma$. In fact, it restricts to the chamber $\mathbb{S}^3\times \DD^2$ preserving boundary and interior. 
In particular, the $4$-dimensional chamber $C^4=(\mathbb{S}^3\times \DD^2)/\S^1_{(1,-1,k,0)}$
associated to the action on $\mathbb{S}^3\times_k \mathbb{S}^2$ is a smooth manifold with unstratified boundary consisting of fixed points of the action, 
and interior with trivial \pig. Notice that $C^4$ has the structure of a $2$-disc bundle over $\mathbb{S}^2$ 
with projection onto the base $\mathbb{S}^2$ induced from that of $\mathbb{S}^3\times \DD^2$ onto its first factor. Such bundles are classified (up to
bundle isomorphism) by $\pi_1(\S^1)\cong \ZZ$, and we may further restrict to non-negative integers by allowing an orientation reversing diffeomorphism
of the base. Finally, they are all topologically distinct and easily identifiable by the fundamental group of their boundary:
$$\partial C =  (\mathbb{S}^3\times \mathbb{S}^1) / \S^1_{(1,-1,k)} \cong \begin{cases}
                                              \mathbb{S}^2 \times \mathbb{S}^1 & k=0 \\
                                              L_k^3 & k\neq 0 
                                             \end{cases}$$
where $L^3_k\cong \mathbb{S}^3/\ZZ_k$ is a Lens space. 
The polar group is $\ZZ_2$ given by a unique reflection, which reconstructs $\Sigma$ as the double of $C$.
\end{ex}

\begin{figure}[h!]
\footnotesize
\begin{center}
\begin{tabular}[c]{cc}
\begin{tikzpicture}[scale=1]
\draw[fill,yshift=0cm]
  (0,0)		circle[radius= 2pt] node[left]{$\SO(3)$} -- node[midway,above] {${12}$} 
  (2,0)		circle[radius= 2pt] node[right]{$\SO(3)$} -- node[midway,right] {${23}$} 
  (1,-1.73) 	circle[radius= 2pt] node[below]{$\SO(3)$} -- node[midway,left] {${13}$} (0,0) ;
\end{tikzpicture}
\end{tabular}
\end{center}
 \caption{ \label{fig Wu} \footnotesize
Coxeter polar data for the $\SO(3)$ action on the Wu manifold. The labeling of the sides corresponds to the block embeddings 
of $O(2)\rightarrow \SO(3)$ in the $ij$ coordinates. }
\end{figure}
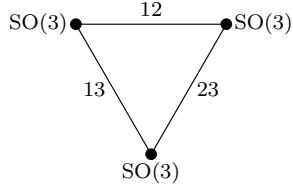

\begin{ex} \label{ex Wu's}
The $\SO(3)$ action from the left on $\SU(3)/\SO(3)=:\W^5$ associated to the symmetric pair $(\SU(3),\SO(3))$ 
is naturally polar with section $\T^2$ and hence Coxeter polar due to \rref{coxeter is natural}. 

Notice that the base point $[e\!\cdot\!\SO(3)] \in \SU(3)/\SO(3)$ 
is a fixed point of the action, where the slice representation is the unique
irreducible $5$-dimensional $\SO(3)$ representation of \eref{ex linear polar}. 
Hence, the \pig~is $\H=\S(O(1)^3)$ and any face isotropy group has to be $O(2)$ since $\K/\H$ is a sphere. 
Furthermore the embedding of $O(2)$ in $\SO(3)$ is determined by an axis, 
when thinking of $\SO(3)$ as acting on $\RR^3$, which is also fixed by the principal isotropy $\H=\S(O(1)^3)$. 
It must hence coincide with one of the coordinate axis, and $O(2)$ be given by one of the three coordinate block embeddings into $\SO(3)$.
Moreover, two consecutive sides of the chamber correspond to different embeddings of $O(2)$ since they need to generate the isotropy group 
at the vertex of intersection. Vertices are thus fixed points of the actions. 
Finally, there are only three vertices, since fixed points $[g\cdot \SO(3)]\in \SU(3)/\SO(3)$ correspond to 
$[g\cdot \SO(3)] \in \N^{\SU(3)}(\SO(3))/\SO(3)\cong \ZZ_3$.
We conclude that the chamber is a triangle 
marked by each of the three block embeddings of $O(2)$ along the sides, as shown in \fref{fig Wu}.

\end{ex}

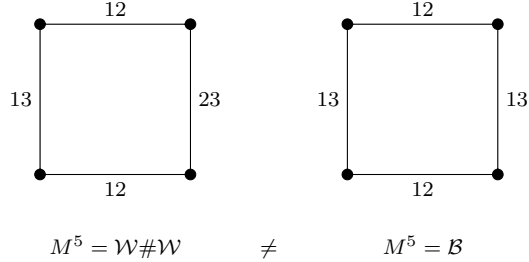
\begin{figure}[h!]
\footnotesize
\begin{center}
\begin{tabular}[c]{ccc}
\begin{tikzpicture}[scale=1]
  \draw[fill,scale=1,xshift=3.5cm,]
  (0,0) 	circle[radius= 2pt] node[above left]{} -- node[midway,above] {${12}$} 
  (2,0)		circle[radius= 2pt] -- node[midway,right] {${23}$} 
  (2,-2)	circle[radius= 2pt] -- node[midway,below] {${12}$} 
  (0,-2) 	circle[radius= 2pt] -- node[midway,left] {${13}$} (0,0)  ;
\end{tikzpicture}
& &
\begin{tikzpicture}
  \draw[fill,scale=1]
  (0,0) 	circle[radius= 2pt] node[above left]{} -- node[midway,above] {${12}$} 
  (2,0)		circle[radius= 2pt] -- node[midway,right] {${13}$} 
  (2,-2)	circle[radius= 2pt] -- node[midway,below] {${12}$} 
  (0,-2) 	circle[radius= 2pt] -- node[midway,left] {${13}$} (0,0)  ;
\end{tikzpicture}
\\
&&\\
$M^5=\W\#\W$& $\neq$ & $M^5=\B$
\end{tabular}
\end{center}
 \caption{ \label{fig Brieskorn Vs Wu} \footnotesize
Coxeter polar data for $\SO(3)$ actions on $M^5$. 
Vertices are fixed points of the action and the labeling of the sides corresponds to the $ij$ coordinate embedding of $O(2)$ in $\SO(3)$.}
\end{figure}

\begin{ex} \label{ex brieskorn}
There is one case among the manifolds listed in \tref{THM classification dim5} whose polarity is only established in an indirect way. This is 
the case of the $\SO(3)$ action on the Brieskorn variety of type $(2,3,3,3)$, which is the zero locus of the polynomial $z_0^2+z_1^3+z_2^3+z_3^3$ 
intersected with the corresponding sphere $\mathbb{S}^7(1)\subset \CC^4$.
This action was constructed in \cite{Hudson} in the classification of $\SO(3)$ actions on $5$-manifolds,  
with prescribed orbit space and isotropy data as in \fref{fig Brieskorn Vs Wu}. 
Using the the Barden-Smale classification of $5$-manifolds the author showed that the manifold
is diffeomorphic to the Brieskorn variety $\B$. Furthermore, it was shown that there is a unique $\SO(3)$ action with such data
up to equivariant diffeomorphism.
In turn, we identify this as compatible Coxeter polar data to which there corresponds a unique Coxeter polar $\SO(3)$-manifold, 
thus establishing the polarity of his construction. Notice that the polar group agrees with polar group at a fixed point and is hence $\D_3$.
$\Sigma$ is a compact surface determined by its Euler characteristic and orientability. 
In fact, it follows from \lref{lem section orientability} that $\Sigma$ is orientable since
any presentation of $\D_3$ with generators given by involutions has even relations only. 
Furthermore, we can regard the quotient map $\Sigma \to C$ under the action of $\Pii$ 
as an orbifold cover (see \cite{Davis}), and compute the Euler 
characteristic of the section from the \textit{orbifold Euler characteristic} \mbox{of $C$} by the formula:
$$\chi(\Sigma) / |\!\Pii\!| = \chi_{orb}(C).$$
In this case $C$ has four points with orbifold group $\D_3$, four sides fixed by an involution and one regular face, 
thus giving 
$\chi_{orb}(C)= 4\cdot 1/6 - 4\cdot 1/2 + 1 = -1/3 $. 
Hence, the section is a bi-torus, $\Sigma=\T^2\#\T^2$.
 
It is interesting to compare the data for this action with the connected fixed point sum of two copies of the Wu manifold $\W\#\W$;
see \fref{fig Brieskorn Vs Wu}. For both $\B$ and $\W\#\W$ the associated section, polar group and chamber coincide, but the action
of $\Pii$ on $\Sigma$ is different. 
\end{ex}

Further fixed point connected sums give rise to polar actions on $\#_{k}\B\#_l\W$. 
The polar group is always the local polar group at a fixed point, $\Pii= \D_3$, hence the section is orientable and corresponds to a connected sum of tori $\#_n \T^2$. 
%
Coxeter polar data can be computed by \lref{fixed point sum}, thus obtaining a polygon marked with \pig~$\H=\S(O(1)^3)$ along the interior, 
vertices given by fixed points of the action, and sides labeled by either of the three coordinate block embeddings of $O(2)$ in $\SO(3)$ without consecutive repetitions.  
Conversely, any such labeling of a polygon gives compatible Coxeter polar data for an $\SO(3)$ action. 
With only two sides this is the linear action on $\mathbb{S}^5$ given by 
the trivial extension of the irreducible $5$-dimensional $\SO(3)$ representation. For three sides, 
the compatibility conditions force all three block 
embeddings to appear once and this corresponds to the isotropy action of $\SO(3)$ on the Wu manifold. With $4$ sides, 
we have two possibilities for the markings, up to conjugation by a single element in $\N(\H)$, 
and these correspond either to the connected fixed point sum $\W\#\W$ or to the Brieskorn variety $\B$, as in \fref{fig Brieskorn Vs Wu}.
One may continue with this kind of reasoning and prove that there is a unique marking with five sides and four possibilities for
the case of six sides.
Finally, an inductive argument on the number of sides shows that the fixed point sums of
$\W$ and $\B$ exhaust all possibilities. 
For this, notice first that if only two block embeddings of $O(2)$ are present in the data, 
then the number of sides must be even and it corresponds to a repeated fixed point connected sum of $\B$ with itself. 
On the other hand, if three different block embeddings appear in the data,
we can find a consecutive sequence of three sides with distinct labels. 
Erasing the middle side of these three and merging two consecutive vertices into one, 
we again have compatible polar data to which we can apply the inductive hypothesis. 
Thus this manifold can be written as a fixed point sum $\#_{k}\B\#_{l}\W$. 
Finally we can perform an additional fixed point sum with the Wu manifold at the vertex where we erased a side, 
which gives back our original data and thus establishes the claim. 
Notice that when we have only one side which has a different block embedding than the rest, we can choose not to remove it but choose a neighboring side
instead, thus keeping the three different block embeddings of $O(2)$ on the resulting data. In this way the induction continues and proves that
the original manifold can be written as a repeated fixed point sum of Wu manifolds only, $\#_n\W$.
We conclude that the $\SO(3)$ actions constructed out of such Coxeter polar data are connected fixed point sums $\#_n \B$ or $\#_n \W$. 
However, while there is a unique action on $\#_n \B$ up to conjugation, this is not the case for $\#_n\W$ since it depends on the possible markings,
thus giving rise to inequivalent actions on the same manifold.

\begin{ex} \label{ex SO3xS1}
The $\SO(3)\times \S^1$ factor-wise action on $\mathbb{S}^3\times \mathbb{S}^2$  
is a product of cohomogeneity one actions. Coxeter polar data is given by a rectangle with four fixed points, 
and sides labeled by singular isotropy groups $\K_{2i+1}= \SO(3)$ and $\K_{2i}= \T^2$.
The \pig~is $\H=\S^1_{1,0}=:\SO(2)$.

We may perform fixed point sums of the latter action with itself, obtaining new polar manifolds.
A generic picture for the Coxeter polar data is given in \fref{fig connected sums SO(3)xS1 on SS3xSS2}. 
Notice that the polar group is always $(\ZZ_2)^2$ and the section $\T^2 \# \cdots \# \T^2$.
\end{ex}
 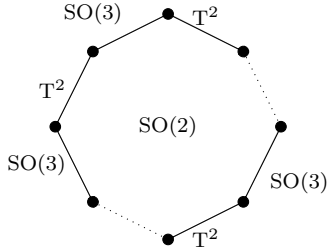
\begin{figure} [h!]
\begin{center}
 \begin{tikzpicture}[scale=1]
 \footnotesize
  \draw[fill]	
  (0,0) circle[radius=2pt] 
	  -- node[above left]{$\SO(3)$} 	(1,1/2) circle[radius=2pt] 		
	  -- node[above]{$\T^2$}		+(1,-1/2) circle[radius=2pt]	
  (0,0)  
	  -- node[left]{$\T^2$} ++(-1/2,-1) circle[radius=2pt] 	
	  -- node[left]{$\SO(3)$} ++(1/2,-1) circle[radius=2pt] 

	  (1,-5/2) circle[radius=2pt] -- node[below]{$\T^2$}	++(1,1/2)circle[radius=2pt] --node[below right]{$\SO(3)$} +(1/2,1) circle[radius=2pt] ;
\draw[dotted]   
  (0,-2) circle[radius=2pt] -- +(1,-1/2) 
  (2,0) circle[radius=2pt]  -- +(1/2,-1) ;
	
\draw (1,-1) node {$\SO(2)$} ;
\end{tikzpicture}
\end{center}
 \caption{ \label{fig connected sums SO(3)xS1 on SS3xSS2}
\footnotesize 
 Polar data for the $\SO(3)\times \S^1$ action on $\#_n \mathbb{S}^3\times \mathbb{S}^2$, where the vertices correspond to fixed points. }
\end{figure}
A diagram with only two vertices is also possible and corresponds to the suspension of the $\SO(3)\times \S^1$ action on $\mathbb{S}^4$, i.e.,
the action on $\mathbb{S}^5$ given by the trivial extension of the product $\SO(3) \times \S^1$ representation on $\RR^3\oplus \RR^2$.

%
%
%
\sect{Classification in low dimensions \label{section Classification in low dimensions}}

We shall consider compact simply-connected smooth manifolds $M^n$, $n\leq 5$, on which a compact connected Lie group $\G$ acts polarly, thus
Coxeter polar, and classify them up to equivariant diffeomorphism. 
The actions are assumed to be almost effective, i.e., with a discrete ineffective kernel. 
As in the previous section, $\H$ will denote the \pig~along a section $\sigma: \Sigma^k \to M^n$ with polar group $\Pii$ and chamber $C$. 

\smallskip
The strategy in order to achieve a classification is to solve the equivalent problem of determining all sets of 
inequivalent Coxeter polar data (see \tref{GZ coxeter reconstruction}). 
To do so, we will consider the possible groups, chambers and compatible group graphs marking the stratifications. We will then
identify these data with the actions in Section \ref{section Examples of polar actions}.

To begin, we notice that the bound on the dimension of the manifold restricts both the cohomogeneity $k$
and the dimension of the group, as follows from the formulas,
\begin{equation} \label{regular dimension count}
\begin{array}{l}
n - k = dim(\G) - dim(\H) \\
\\
dim(\H) \leq dim(O(n-k)).
\end{array}
\end{equation}
The first is a simple dimension count on the decomposition $T_pM = T_p(\G\cdot p) \oplus S_p$ at a regular point. 
The second comes from the fact that the representation of the principal isotropy group in the orbit direction 
corresponds to the 
restriction of the adjoint representation of $\G$ to $\H$, and its kernel to the ineffective kernel of the $\G$-action. 
Therefore, we have $\H \subset O(n-k)$.

Moreover, any compact connected Lie group $\G$ can be taken, up to a finite cover, to be of the form
\begin{equation} \label{group G up to covers}
\G = \T^l\times \G_1 \times \cdots \times \G_k 
\end{equation}
for a torus $\T^l$ and simply-connected simple factors $\G_i$. 
Together with the bound on the dimension, this yields finitely many possibilities for the group.

Notice that if the action does not have singular orbits, the section will agree with the chamber, $C=\Sigma$, which is thus a smooth manifold without boundary 
constantly marked by a principal isotropy group $H$. This is an example of an imprimitive action for which the reconstructed 
manifold is a product of the section with a regular orbit, $M= \G/\H \times \Sigma$. These Coxeter polar actions will be called \textit{trivial}.

Since $M$ is simply-connected, the orbit space $M/\G$ is compact and simply-connected as well.
Furthermore, it is isometric to a chamber $C^k\subset \Sigma$ 
which has the structure of a smooth Coxeter orbifold, i.e.,  
one modeled on the quotient of euclidean $k$-space by finite Coxeter groups. 
In particular, there will be no singularities along the interior and we may ``round the corners'' along the boundary to obtain
a homeomorphism of $C$ with a topological manifold with boundary. 

\begin{prop} \label{2,3 quotients}
A compact simply-connected $n$-manifold with boundary is homeomorphic to $\DD^2$ if $n\!=\!2$, 
or a sphere with finitely many discs removed $\mathbb{S}^3 \backslash \sqcup_{i=0}^r \DD^3_i$, $r \in \NN$, if $n\!=\!3$.
\end{prop}
The first case follows from the classification of closed surfaces, when capping the circles in the boundary with closed discs.
The case $n=3$ is a consequence of Poincar\'e duality, as shown in \cite{Simas}. 

Once the orbit space is determined as a topological manifold with boundary, 
we must still address its stratification and what the 
markings of the strata by isotropy subgroups could be.
This is particularly simple in the case of polar actions of maximal cohomogeneity, i.e., those with codimension one section.

\begin{prop} \label{cohomogeneity n-1}
A polar action on a closed $n$-manifold of cohomogeneity \mbox{$n\!-\!1$} is a circle action with 
orbit space given by a smooth manifold with boundary marked as fixed points. 
Conversely, to each smooth manifold with boundary corresponds a unique polar circle action on a closed manifold.

The associated manifold is simply-connected if and only if the orbit space is also simply-connected. 
\end{prop}
\begin{pf}
 If the cohomogeneity of the action is $n-1$ we have by \eqref{regular dimension count} that
 the principal isotropy has to be discrete, since $O(1)$ is. and the group $\G$ one-dimensional.
 When made effective this is an $\S^1$-action with trivial \pig~and $\K=\S^1$ as 
 the only other possible singular isotropy subgroup. 
 This corresponds to fixed points on $M$, which constitute the boundary of the chamber. 
 
As a converse, any (simply-connected) smooth manifold with boundary $C^{n-1}$ 
gives rise to a unique Coxeter polar $\S^1$-manifold $M^n=M(C,\S^1(C))$ by \tref{GZ coxeter reconstruction}. 
In particular, if this boundary were empty we are in the trivial case $M= M' \times \S^1$, which is not simply connected. 
Furthermore, if the boundary is not trivial, Seifert-Van Kampen's theorem easily 
implies that the fundamental group of the reconstructed manifold is that of the chamber $C$.
Indeed, notice that the interior of $C$  
gets reconstructed as $\text{int}(C)\times \S^1$, while an $\epsilon$-collar neighborhood of a boundary component of $C$
gives rise to $\partial_0C\times \DD^2\cong D_{\epsilon}(\partial_0C) \subset M$, from the uniqueness of the reconstruction process. 
Each such piece is to be glued with the interior along the common intersection 
$(0,\epsilon) \times \partial_0 C \times  \S^1$, which cancels the fundamental group of the $\S^1$ factor.
%
We then continue by induction on the boundary components. 
\qed
\end{pf}

\begin{rem*}
A polar circle action is \textit{fixed point homogeneous}, i.e., the group $\S^1$ acts transitively on the sphere orthogonal to the fixed point stratum.

The polar group is $\ZZ_2$, generated by the unique reflection in $\S^1$ 
and the section is the double of the chamber, $\Sigma =\mathcal{D}(C)$. 
\end{rem*}

\begin{cor} \label{classification dim3}
There is only one effective compact simply-connected cohomogeneity two polar $3$-manifold, 
given by a linear $\S^1$-action on $\mathbb{S}^3$. 
\end{cor}
%
We now describe the classification in dimension $4$.
\begin{thm} \label{classification dim4}
Let $G$ be a compact connected Lie group, with a non-trivial and and almost-effective polar action on a compact simply-connected $4$-manifold $M$ 
of cohomogeneity $k\geq 2$. Then $M$ is equivariantly diffeomorphic to one of the following:
\begin{enumerate}
\item[a)] $\T^2$-actions as described in \eref{Torus cohom2}; 
\item[b)] The linear $\S^1$ and $\SO(3)$ actions on $\mathbb{S}^4$ as in \eref{ex linear polar}; 
\item[c)] The linear circle action on 
the first factor of $\mathbb{S}^2 \times \mathbb{S}^2$, or their connected fixed-point sums $\#_n\mathbb{S}^2\times \mathbb{S}^2$.
  \end{enumerate}
\end{thm}
\begin{pf}
The circle actions corresponds to the case $k=3$ which was already considered in \pref{cohomogeneity n-1}.
 Additionally we identify the reconstructed manifolds to be $\mathbb{S}^4$ when $C=\DD^3$, $\mathbb{S}^2\times \mathbb{S}^2$ for $C= \DD^1 \times \mathbb{S}^2$,
 and fixed point sums of the latter in the remaining cases. In fact, 
 each fixed point connected sum of the previous $\mathbb{S}^2\times \mathbb{S}^2$ adds one sphere component 
 to the boundary of the chamber (see \pref{fixed point sum}), allowing to retrieve any 
 possible chamber $C= \mathbb{S}^3 \backslash \sqcup_{i\leq N} (\DD^3)_i$. 
 
If the cohomogeneity is $2$, either the \pig~is discrete and hence $dim(\G)=2$, which implies $\G=\T^2$,
or we have a one dimensional \pig~and $\G=\SU(2)$. 
For this last case, since 
there are no subgroups of $\SU(2)$ of intermediate dimension, the chamber has to be a disc $\DD^2$ 
with un-stratified boundary marked by $K=\SU(2)$.
As one can check, this coincides with the polar data for the effective $\SO(3)$ action on $\mathbb{S}^4$ 
given by its irreducible $3$-dimensional representation plus two fixed directions.
\qed
\end{pf}

\begin{cor}
 A closed simply-connected $4$-manifold admitting a polar action is diffeomorphic to $\mathbb{S}^4$ or 
 a connected sum of copies of $\pm \CC\PP^2$ and $\mathbb{S}^2\times \mathbb{S}^2$.
\end{cor}
The previous corollary combines 
our classification 
and the work on cohomogeneity one manifolds in \cite{Hoelscher}. It should be mentioned that this is true even for non-polar actions 
as observed in \cite{GeRad}.

We now address the classification of polar actions by non-abelian groups in \mbox{dimension $5$}, starting with those
by $\SU(2)$ and $\SO(3)$.

\begin{thm} \label{thm non abelian 5mflds}
A non-trivial polar $\SU(2)$ action on a closed simply-connected $5$-manifold is equivariantly diffeomorphic to one of the following:
\begin{enumerate}
 \item The linear polar actions on $\mathbb{S}^5$ listed in \eref{ex linear polar}; 
 \item The actions on $\mathbb{S}^3\times_k \mathbb{S}^2$ described in \eref{ex SU(2) of main example};
 \item The $\SO(3)$-action on the Brieskorn variety $\B$ of type $(2,3,3,3)$, the left $\SO(3)$-action on the Wu manifold $\W=\SU(3)/\SO(3)$, 
 or connected fixed point sums $\#_n\B$ and $\#_m\W$;
 \item The linear $\SO(3)$ action on the first factor of $\mathbb{S}^3\times \mathbb{S}^2$ and repeated connected fixed point sums, 
 ${\#}_n \mathbb{S}^3\times \mathbb{S}^2$. 
\end{enumerate}
\end{thm}
\begin{pf}
We distinguish between two main cases depending on whether the \pig~is discrete or one dimensional. In the latter case, we have a cohomogeneity three action, 
with singular isotropy given by the whole group since there are no subgroups of $\SU(2)$ of intermediate dimension. It follows that $\H=\S^1$, and
that the chamber is a smooth manifold with unstratified boundary marked as fixed points. 
A basic example is the $\SO(3)$ linear action on $\mathbb{S}^5$ given by the trivial
extension of its $3$-dimensional representation, which corresponds to chamber $C=\DD^3$. The remaining cases correspond to the action
on the first factor of $\mathbb{S}^3\times \mathbb{S}^2$ with chamber $\DD^1\times \mathbb{S}^2$, or their connected fixed point sums which have 
$C=\mathbb{S}^3\backslash \bigsqcup_i \epsilon \DD^3_i$.

We now assume that $\H$ is discrete. Hence the action has cohomogeneity two and $C=\DD^2$. 
The possible singular isotropy groups are $\K=\S^1$, $\Pin(2)$ or $\SU(2)$.
If the boundary is unstratified and constantly marked by $\K=\S^1$ or $\Pin(2)$, this gives imprimitive actions 
which reconstruct $\S^3$ as $\K$-manifolds, while as $\SU(2)$-manifolds they are equivariantly diffeomorphic to
$$	M=\SU(2)\times_{\K}\mathbb{S}^3	.	$$
Since $\K$ acts freely on the product, we conclude from the associated long homotopy sequence 
that, if $\K$ is not connected, then its quotient $M$ is not simply-connected.
If we assume $\K=\S^1$ and $\H=\ZZ_k$, this is uniquely identified with the examples treated in \ref{ex SU(2) of main example},
which correspond to actions on the biquotient $\mathbb{S}^3\times_k\mathbb{S}^2$.
If on the other hand, the whole boundary corresponds to fixed points of the $\SU(2)$ action, the \pig~is then trivial and we recognize this as
the linear $4$-dimensional $\SU(2)$-representation extended to act on $\mathbb{S}^5$.

For the remaining cases the chamber has stratified boundary, with isolated vertices corresponding to fixed points of the action, 
in order to satisfy \eqref{coxeter rigidity2}.
The slice representation at such points must be the irreducible $5$-dimensional one of \eref{ex linear polar}.2. 
These cases correspond to the action on $\mathbb{S}^5$ by the (trivial extension of the) $\SO(3)$ \mbox{$5$-dimensional} representation, 
the isotropy action of $\SO(3)$ on the Wu manifold, the $\SO(3)$ action on the Brieskorn variety, and their connected fixed point sums. 
We refer to the discussion in Examples \ref{ex Wu's}, \ref{ex brieskorn} and to Figures \ref{fig Wu} and \ref{fig Brieskorn Vs Wu}.
\qed
\end{pf}

\begin{rem}
General $\SO(3)$ actions on $5$-manifolds were previously studied in \cite{Hudson} and more recently in \cite{Simas} who corrected 
and extended the classification to that of (almost effective) $\SU(2)$ actions.
An alternative proof of \tref{thm non abelian 5mflds} follows from Theorem D in \cite{Simas}, after identifying the polar examples 
and discarding the non-polar ones. The latter correspond to the linear diagonal $\SO(3)$ action on $\mathbb{S}^5$ 
and to the left $\SU(2)$-action on the Wu manifold.
The $\SO(3)$ action on $\mathbb{S}^5\subset \RR^{3+3}$ has orbit space a disk with
boundary marked by $\S^1$ isotropy, and trivial principal isotropy along the interior. 
This gives compatible imprimitive Coxeter polar data for an $\SO(3)$ action, but corresponds to the (clearly polar) $\SO(3)$ action on 
$\mathbb{S}^3\times_{k=2}\mathbb{S}^2\cong \mathbb{S}^3\times\mathbb{S}^2$, presented in \eref{ex SU(2) of main example} as an ineffective $\SU(2)$ action.
In the same way, the left $\SU(2)$-action on $\W$ has orbit space a disc with constant $\S^1$ isotropy along the bounding circle and
trivial isotropy in the interior, which gives rise to Coxeter polar data reconstructing $\mathbb{S}^3\times \mathbb{S}^2$ instead of $\mathbb{S}^5$.
\end{rem}

For general non-abelian group actions we have the following, which thus finishes the proof 
of \tref{THM classification dim5} in the introduction.  
\begin{thm} \label{thm classification dim5}
Let $\G$ be a compact connected non-abelian Lie group, 
with a non-trivial and almost-effective polar action on a compact simply-connected $5$-manifold $M$ of
cohomogeneity $k=2,3$. Then $M$ is equivariantly diffeomorphic to one of the following:
\begin{enumerate}
\item[a)] The $\SU(2)$ polar actions listed in \tref{thm non abelian 5mflds};
\item[b)] The family of $\SU(2) \times \S^1$-actions on $\mathbb{S}^3 \times_k \mathbb{S}^2$ as in \eref{ex SU(2)xS1 of main example above}; 
\item[c)] The linear factor-wise $\SO(3)\!\times\! \S^1$ action on $\mathbb{S}^3 \times \mathbb{S}^2$, 
and their connected fixed-point sums;
\item[d)] Linear polar actions by $\U(2)$ and $\SO(4)$ on $\mathbb{S}^5$.
\end{enumerate}
\end{thm}

\begin{pf}
From the same considerations as in the classification in dimension four, we have that a cohomogeneity $k=3$ action corresponds
to an $\SU(2)$ action with one dimensional \pig.  

If $k=2$ we have $C\cong \DD^2$ and $\H\subseteq O(3)$ up to a finite cover, leaving the following possibilities: 
\begin{equation}\label{group and pig}
(\G,\H_0)=(\SU(2),e),(\SU(2)\!\times\! \SU(2),\Delta\! \SU(2)),(\SU(2)\!\times\!\S^1,\S^1).
\end{equation}
Here, besides listing the possible groups by the dimension bound \eqref{regular dimension count} for a given $\H_0\subset \G$,
we have discarded $\G=\SU(2)\times \T^3$ with $\H_0=\SU(2)$,
since it would not be almost-effective.
The $\SU(2)$ actions of cohomogeneity $2$ and $3$ have been addressed in 
\tref{thm non abelian 5mflds}.

For $(\G,\H_0)=(\SU(2)\times \SU(2),\SU(2))$, notice that the effectiveness of the $\G$-action implies that $\H_0$ is embedded diagonally. 
Furthermore, there are no subgroups of intermediate dimension between $\G$ and $\H_0$, 
implying that the boundary of the orbit space consists of fixed points by $\G$. 
In particular, we have that $\H$ must be connected, since $\G$ corresponds to a face isotropy group and satisfies $\G/\H=\mathbb{S}^3$. 
Thus the polar data is a disc with interior regular isotropy $\Delta \SU(2)$ 
and boundary consisting of fixed points.
If we pass to the quotient under the normal subgroup $\ZZ_2\cong \{\pm 1\} \subset \Delta \SU(2)$, we have an effective $\G=\SO(4)$ 
action with \pig~$\H=\SO(3)$. This equivalent effective data is recognized as the linear $\SO(4)\subset \SO(6)$ action on $\mathbb{S}^5$.

If $(\G,\H_0)=(\SU(2)\!\times\!\S^1,\S^1)$, we start by determining the ways in which $\H_0$ may lie in $\G$. 
This cannot be given by the inclusion of the second factor, as it would not be almost effective. 
It must hence be 
$\H_0=\S^1_{(1,k)}\subset \SU(2)\times \S^1$, 
for $k \in \ZZ$, since a general circle $\S^1_{(q,k)}$ 
would have normal ineffective kernel $\ZZ_q\cong <e^{i2\pi/q}>\subset \S^1$ which
can be quotiented out and still retrieve an $\SU(2)\times \S^1$ action.
We distinguish between the inclusion of the first factor only ($k=0$) and the diagonal case ($k\neq 0$).

In the case of $\H_0=\S^1_{(1,k)}$, $k\neq 0$, the singular isotropy groups can be $\K\cong \T^2$, $\Pin(2)\times \S^1$,
or the whole group $\SU(2)\!\times\!\S^1$. We claim that, in any case, the boundary will have constant isotropy. 
Notice that there is only one choice for the maximal torus inside the different possible singular isotropy subgroups of a face, 
namely, the one given by the product of the circle subgroups obtained from factor-wise projection of $\H_0\cong \S^1\subset \SU(2)\!\times\!\S^1$.
This rules out the possibility of a further stratification of the boundary, since consecutive sides
have to be labeled by different isotropy groups and generate the isotropy group of the vertex in common. 
In the case of unstratified boundary constantly marked by $\K=\T^2$ or $\Pin(2)\times \S^1$ we have imprimitive actions which reconstruct $\mathbb{S}^3$ 
as a $\K$-manifold, while as a $\SU(2)\times \S^1$-manifold they give
$$(\SU(2)\times \S^1) \times_{\K} \mathbb{S}^3.$$
In particular, $\K$ has to be connected in order for $M$ to be simply-connected.
Thus we can assume $\K=\T^2$ which is the case discussed in \eref{ex SU(2)xS1 of main example above}.
The remaining case with unstratified boundary equal to the fixed point set can only happen if the \pig~is $\H= \S^1_{(1,1)}$,
in order to satisfy $\K/\H\cong \mathbb{S}^3$. This is identified as the linear effective action on $\mathbb{S}^5$ of the group
$\U(2)\cong \SU(2)\cdot \S^1$ acting as a subgroup of $\SO(4)\subset \SO(6)$. 

Finally, we turn to the case where the principal isotropy is $\H_0=\S^1_{(1,0)}$, which corresponds to an effective $\SO(3)\times \S^1$ action 
(since we can quotient out the normal subgroup $\{(\pm 1,1) \} \subset \H_0 \subset \SU(2)\times \S^1$).
If the boundary is unstratified it is thus constantly marked by either $\K=\T^2$, $\SO(3)$ or $\SO(3)\times \S^1$. 
The case with $\K=\SO(3)$ gives an imprimitive action that reconstructs a non-simply connected manifold.
The other two can be treated as before,
and are easily identified as the action on $\mathbb{S}^3\times_{k=0} \mathbb{S}^2=\mathbb{S}^3\times \mathbb{S}^2$ in \eref{ex SU(2)xS1 of main example above}, 
and the linear action on $\mathbb{S}^5$ given by the trivial extension of the $\SO(3)\times \S^1$ representation in $\RR^3\oplus \RR^2$. 
Finally, if the boundary is stratified the marking of the sides must alternate between $\K=\T^2$ and $\SO(3)$. 
These actions were the ones discussed in \eref{ex SO3xS1}.
\qed
\end{pf}

%
\sect{Torus actions \label{section Polar $T^2$ actions on $M^5$}}
In this section we treat in detail the case of polar effective $\T^2$-actions on compact simply-connected $5$-manifolds. 
%
%
We will prove \tref{THM T2 in M5} 
by determining the admissible Coxeter polar data, addressing the possible isotropy groups, orbit spaces and marked stratifications.
In a second step we identify the actions concretely, 
which requires the description of a few basic examples and the combination of them by means of equivariant surgeries.

First, observe that the possible isotropy types $(\G_x)$ are known and correspond to \mbox{$\{e\}=\H \subseteq \S^1 \cong \K \subseteq \T^2$}, 
where $\S^1$ ranges over all the possible circle subgroups inside the torus, and we omit parenthesis as conjugation is trivial.
In fact, the effectiveness of the action implies that the \pig~is trivial since the group is abelian. 
A next to minimal isotropy subgroup has to be a circle,
since $\K/\H$ is a sphere. 
Moreover, since $\K_0$ and $\H$ generate $\K$, further isotropy groups will have to be connected leaving as only possibility $\K=\G=\T^2$, i.e.,
fixed points of the action.

Secondly, notice that at a fixed point the $5$-dimensional slice representation coincides with
the standard one by the maximal torus $\T^2 \subset \SO(4) \subset \SO(5)$. 
In particular, the existence of a fixed direction implies that fixed points are not isolated, but form closed geodesics. 
For a point with $\S^1$ isotropy, the orbit is a circle and the slice must be $4$-dimensional. The slice representation is the usual $2$-dimensional one 
by $\S^1$ extended trivially.

The possible orbit spaces are topological $3$-manifolds with boundary, already characterized in \pref{2,3 quotients} 
as the complement of a finite number of disjoint $3$-discs in a $3$-sphere,
$$C=M/\G =\mathbb{S}^3 \backslash (\sqcup_N~ \DD^3).$$

Each component of the boundary $\mathbb{S}^2 \subseteq \partial (M/G)$ will be the union of one or more faces with $\S^1$-isotropy, 
separated by closed geodesics consisting of fixed points. 
We can encode the boundary stratification and marking information in an unoriented graph, 
defining vertices to be the $2$-dimensional faces and edges to be the closed geodesics fixed by the action.
It is well defined, since these fixed geodesics do not intersect each other, leaving only two adjacent faces corresponding to their endpoints. 
Also, as any simple closed curve separates the \mbox{$2$-sphere} into two pieces, this graph is in fact a tree. 
The marking of the chamber labels each vertex with a pair of coprime integers $(m_i,n_i) \in \ZZ^2$, 
indicating the slope of the corresponding isotropy circle subgroup.
The compatibility condition \eqref{coxeter rigidity2} corresponds to consecutive vertices 
having isotropy subgroups with trivial intersection, or equivalently that 
$$ \left|\!\begin{array}{cc}
m_i & n_i \\
m_{i+1} & n_{i+1} 
 \end{array}\!\right|=\pm 1 .$$
We call this graph the \textit{boundary tree}, which is a more convenient way of describing the Coxeter polar data, adapted to this case. 
In particular, any finite tree with a compatible marking of slopes and possibly many connected components, 
will give rise to a unique connected Coxeter polar $\T^2$-manifold of dimension five. 
This reconstructed manifold need not be simply-connected, a problem we shall address later. 

We first focus on actions having only one boundary component for the orbit space, i.e., a connected boundary tree. 
The simplest case is the one without fixed points, and hence only one face marked by a circle isotropy subgroup $\K$.
It corresponds to a tree with only one vertex.
This action is imprimitive and reconstructs $\mathbb{S}^4$ when restricted to $K=\S^1$ with the action by rotations 
in $\RR^2$, trivially extended to $\RR^5$. As a $\T^2$ polar manifold it yields 
$$\mathbb{S}^4 \times_{\K} \T^2 = \mathbb{S}^4 \times_{\K} (\K\times \mathbb{S}^1)= \mathbb{S}^4\times \mathbb{S}^1.$$ 
Although this manifold is not simply-connected it will be one of the building blocks in our construction.

A boundary tree with two vertices connected by a single edge corresponds to a pair of neighboring faces separated by one geodesic 
consisting of fixed points. 
We can assume that the circle isotropy subgroups marking the faces 
are the canonical ones of slopes $(1,0)$ and $(0,1)$, up to a $\T^2$ automorphism. 
The manifold is recognized as the linear $\T^2$ action on $\mathbb{S}^5$ given by the standard embeddings $\T^2\subset \SO(4)\subset \SO(6)$.

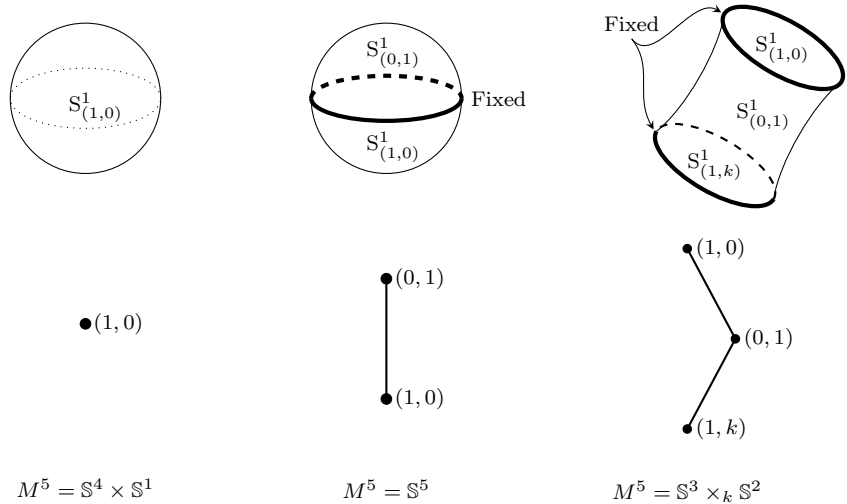
\begin{figure}[h!]
\begin{center}
\newcommand{\length}{3} 
\newcommand{\radius}{1} 
\newcommand{\altura}{2}
\newcommand{\ecc}{0.3cm }

\footnotesize 

\begin{tikzpicture}
\draw[]	(0,0) circle[radius=1] ;
 \draw (0.14,-0.12) node {$\S^1_{(1,0)}$};
 \draw[thin,dotted] ellipse(1cm and 0.4cm); 
\begin{scope}[xshift=4cm]
\draw (0,0) circle (1cm);
\draw[ultra thick] (-1,0) arc (180:360:1cm and \ecc) node[right] {Fixed} ;
\draw (0.1,-0.65) node[] {$\S^1_{(1,0)}$};
\draw[ultra thick, dashed] (-1,0) arc (180:0:1cm and \ecc); 
\draw (0.1,0.6) node[] {$\S^1_{(0,1)}$};
\end{scope}
\begin{scope}[xshift=10.5cm,yshift=1cm]
\begin{scope}[rotate=-30,scale=0.9,yshift=-3cm]
 \draw[ultra thick] (0,\altura) arc(0:360:1cm and 0.4cm);
\draw[dashed,thick] (0,0) arc(0:180:1cm and 0.4cm);
\draw[ultra thick] (0,0) arc(0:-180:1cm and 0.4cm);	
\draw (0,0) .. controls +(-0.2,0.5) and +(-0.2,-0.5) .. (0,\altura)
      (-2,0) .. controls +(0.2,0.5) and +(0.2,-0.5) .. (-2,\altura);
\draw node at (-1,0) {$\S^1_{(1,k)}$} ;
\draw node[right] at(-1.2,0.37*\altura) {$\S^1_{(0,1)}$}; 
\draw node at(-1,\altura){$\S^1_{(1,0)}$};
\newcommand{\pos}{(-3,1)}
    \draw[color = black] \pos node[above]{Fixed};
    \draw[-stealth] \pos .. controls +(0.5,-0.2) and +(-0.5,+0.3)	.. (-2.05,0);
    \draw[-stealth] \pos .. controls +(0.5,-0.2) and +(-0.5,0) 	.. (-2.05,\altura);
\end{scope}
\end{scope}

\begin{scope}[yshift=-2cm]
 
\draw (0,-3.2) node {$M^5=\mathbb{S}^4\times \mathbb{S}^1$} ; 
\draw[white,yshift=-1cm] (0,1)--(0,0)--(0,-1) 
			(1,0)--(0,0)--(-1,0) ;
\filldraw[yshift=-1cm] (0,0) circle[radius=2pt] node[right] {$(1,0)$} ; 

\begin{scope}[scale=1,rotate=0,xshift=4cm,yshift=-0.4cm]
 \draw[thick] (0,0) node[right] {$(0,1)$} -- (0,-1.6) node[right] {$(1,0)$};
\filldraw (0,0) circle[radius=2pt] (0,-1.6) circle[radius=2pt] ;
\draw (0,-2.8) node {$M^5=\mathbb{S}^5$} ; 
\end{scope}

 \begin{scope}[scale=0.8,rotate=0,xshift=10cm,yshift=0cm]
\draw[thick] (0,0) node[right]{$(1,0)$} -- +(0.8,-1.5) node[right]{$(0,1)$} -- +(0,-3) node[right]{$(1,k)$};
\filldraw (0,0) circle[radius=2pt] +(0.8,-1.5) circle[radius=2pt] +(0,-3)circle[radius=2pt] ;
\draw (0,-4.04) node {$M^5=\mathbb{S}^3\times_k \mathbb{S}^2$} ; 
\end{scope}
\end{scope}
\end{tikzpicture}
\end{center}
\caption{ \footnotesize \label{fig basic boundary trees}
Coxeter polar data and their corresponding boundary trees for the basic $\T^2$ actions 
on $5$-manifolds with orbit space homeomorphic to the $3$-disc.}
\end{figure}

The $\T^2$ action on $\mathbb{S}^3\times_k\mathbb{S}^2=\mathbb{S}^3\times \mathbb{S}^3/\S^1_{(1,1,k,0)}$ as described in \eref{ex T3 case main example}
have boundary tree given by a linear sequence of three vertices labeled as in \fref{fig basic boundary trees}.
Conversely, any $\T^2$ Coxeter polar action with the same underlying boundary tree can be 
normalized to have the same marking by isotropy subgroups, as follows. 
First pick two consecutive vertices and assume that their isotropy groups are the circles of slopes
$(1,0)$ and $(0,1)$, up to a $\T^2$ automorphism. 
Furthermore, we can choose the middle vertex to be labeled $(0,1)$, so that the remaining one will have isotropy subgroup of the form $(1,k)$, 
for some integer $k\in \ZZ$, since it is a singular isotropy subgroup with trivial intersection with the circle $(0,1)$.  

These manifolds can be combined by taking fixed point connected sums to obtain new $\T^2$ polar manifolds.
The resulting action may vary depending on the point where it is performed, as opposed to the usual topological connected sum. 
Once compatible points are specified the gluing is determined by the marking of isotropy groups, 
since orientations for the manifolds have to be chosen so that both have outward (or inward) pointing normal vector field along
the maximal faces of the boundary. We illustrate these subtleties in \fref{fig surgeries}, where  
two almost identical pairs of manifolds are glued, and in fact only one summand differs from its analogous by a $\T^2$ automorphism. 

\newcommand{\length}{3} 
\newcommand{\radius}{1} 
\newcommand{\altura}{2}
\begin{figure}[p]	

\begin{center}
\begin{tabular}{cc}
 \begin{tikzpicture}[yshift=0cm,xshift=-2cm,scale=0.88]
 \draw (3,1) node {$\neq$};
 \begin{scope}[xshift=-2cm]
\begin{scope}[rotate=-30]					
 \draw[ultra thick] (0,\altura) arc(0:360:1cm and 0.4cm);
\draw[dashed,thick] (0,0) arc(0:180:1cm and 0.4cm);
\draw[ultra thick] (0,0) arc(0:-180:1cm and 0.4cm);	
\draw (0,0) .. controls +(-0.2,0.5) and +(-0.2,-0.5) .. (0,\altura)
      (-2,0) .. controls +(0.2,0.5) and +(0.2,-0.5) .. (-2,\altura);

\draw node at (-1,0) {$\S^1_{(1,k)}$} ;
\draw node[right] at(-1.2,0.37*\altura) {$\S^1_{(0,1)}$}; 
\draw node at(-1,\altura){$\S^1_{(1,0)}$};

\draw[dotted,fill=white] (0,\altura)  circle[radius=4pt] node[right] {\,\#};
\end{scope}    
 \begin{scope}[xshift=3.48cm,yshift=0.74cm,rotate=-30]		
  \draw[ultra thick] (0,\altura) arc(0:360:1cm and 0.4cm);
\draw[dashed,thick] (0,0) arc(0:180:1cm and 0.4cm);
\draw[ultra thick] (0,0) arc(0:-180:1cm and 0.4cm);	
\draw (0,0) .. controls +(-0.2,0.5) and +(-0.2,-0.5) .. (0,\altura)
      (-2,0) .. controls +(0.2,0.5) and +(0.2,-0.5) .. (-2,\altura);
\draw node at (-1,0) {$\S^1_{(0,1)}$} ;
\draw node[right] at(-1.2,0.4*\altura) {$\S^1_{(1,0)}$}; 
\draw node at(-1,\altura){$\S^1_{(q,1)}$};
\draw[dotted,fill=white] (-2,0)  circle[radius=4pt]; 
\end{scope}
 \end{scope}
\end{tikzpicture}

& 

\begin{tikzpicture}[yshift=0cm,scale=0.88]

\begin{scope}[rotate=-30]
 \draw[ultra thick] (0,\altura) arc(0:360:1cm and 0.4cm);
\draw[dashed,thick] (0,0) arc(0:180:1cm and 0.4cm);
\draw[ultra thick] (0,0) arc(0:-180:1cm and 0.4cm);	
\draw (0,0) .. controls +(-0.2,0.5) and +(-0.2,-0.5) .. (0,\altura)
      (-2,0) .. controls +(0.2,0.5) and +(0.2,-0.5) .. (-2,\altura);

\draw node at 		(-1,0) {$\S^1_{(1,k)}$} ;
\draw node[right] at	(-1.2,0.37*\altura) {$\S^1_{(0,1)}$}; 
\draw node at		(-1,\altura){$\S^1_{(1,0)}$};
\draw[dotted,fill=white] (0,\altura)  circle[radius=4pt] node[right] {\,\#};
\end{scope}
\begin{scope}[xscale=-1,xshift=-1.7cm,yshift=-0.29cm]
 \draw[ultra thick] (0,\altura) arc(0:360:1cm and 0.4cm);
\draw[dashed,thick] (0,0) arc(0:180:1cm and 0.4cm);
\draw[ultra thick] (0,0) arc(0:-180:1cm and 0.4cm);	
\draw (0,0) .. controls +(-0.2,0.5) and +(-0.2,-0.5) .. (0,\altura)
      (-2,0) .. controls +(0.2,0.5) and +(0.2,-0.5) .. (-2,\altura);

\draw node at (-1,0) {$\S^1_{(1,q)}$} ;
\draw node[right] at(-0.7,0.45*\altura) {$\S^1_{(0,1)}$}; 
\draw node at(-1,\altura){$\S^1_{(1,0)}$};

\draw[dotted,fill=white] (0,\altura)  circle[radius=4pt] ;
\end{scope}
\end{tikzpicture}
\\ 

\vspace{0.4cm} & 
\\

\begin{tikzpicture}[scale=0.88]
 \begin{scope}[xshift=0.5cm] 
\foreach \n in {0,1} 
\draw[dashed,thick] (0,\n*\altura) arc(0:180:1cm and 0.4cm);
\draw[ultra thick] (0,2*\altura) arc(0:180:1cm and 0.4cm);	
\foreach \n in {0,1,2} 
\draw[ultra thick] (0,\n*\altura) arc(0:-180:1cm and 0.4cm);
\foreach \n in {0,1} 
\draw (0,\n*\altura) .. controls +(-0.2,0.5) and +(-0.2,-0.5) .. (0,{(1+ \n)*\altura})
      (-2,\n*\altura) .. controls +(0.2,0.5) and +(0.2,-0.5) .. (-2,{(1+ \n)*\altura});

\draw node at (-1,0) {$\S^1_{(1,k)}$} ;
\draw node[right] at(0.5,0.5*\altura) {$\S^1_{(0,1)}$}; 
\draw node[right] at(0.5,1.5*\altura) {$\S^1_{(1,0)}$};
\draw node at(-1,2*\altura){$\S^1_{(q,1)}$};
\draw[-stealth] (0.5,0.5*\altura) -- (- 0.5,0.5*\altura) ;
\draw[-stealth] (0.5,1.5*\altura) -- (- 0.5,1.5*\altura) ;
 \end{scope}
 \end{tikzpicture}
 
 & 
 
 \begin{tikzpicture}[scale=0.88]
\begin{scope}[yscale=-1.1,xshift=0.5cm]
\renewcommand{\altura}{2}
\draw[ultra thick]	(1,-\altura)  arc(0:360:1cm and 0.4cm)	;

\draw[ultra thick] (-0.5,0) arc(0:180:1cm and 0.4cm) ;
\draw[dashed, thick]	(-0.5,0) arc(0:-180:1cm and 0.4cm)	;

\draw[ultra thick,xscale=-1] (-0.5,0) arc(0:180:1cm and 0.4cm) ;
\draw[dashed, thick,xscale=-1]	(-0.5,0) arc(0:-180:1cm and 0.4cm)	;
\draw (-0.5,0) .. controls +(0,-1.2) and +(0,-1.2) ..(0.5,0);
\draw (-2.5,0) .. controls +(0,-1) and +(0,1)..(-1,-\altura);
\draw (2.5,0)  .. controls +(0,-1) and +(0,1)..(1,-\altura);

\draw (-1.5,0) node {$\S^1_{(1,k)}$}
(1.5,0) node {$\S^1_{(1,q)}$}
(0,-\altura) node {$\S^1_{(1,0)}$}
(0,-\altura +0.8) node {$\S^1_{(0,1)}$};
\end{scope}
\end{tikzpicture}

 \\ 
 
 \vspace{0.4cm} & 
\\

 \begin{tikzpicture}
\begin{scope}[xshift=-1cm,yshift=1cm,scale=0.8]
\draw[thick] (0,0) -- (0.8,1.5) node[right]{$(q,1)$};
\draw[thick] (0,0) node[left]{$(1,0)$} -- +(0.8,-1.5) node[right]{$(0,1)$} -- +(0,-3) node[left]{$(1,k)$};
\filldraw (0,0) circle[radius=2pt] +(0.8,-1.5) circle[radius=2pt] +(0,-3)circle[radius=2pt] +(0.8,1.5) circle[radius=2pt];
\end{scope}
\end{tikzpicture}

& 
\begin{tikzpicture}
\begin{scope}[xshift=0.5cm,yshift=-8cm,scale=1]
\draw[white] (0,2)--(0,-2);
\draw[thick] (0,1.2) node[left]{$(1,0)$} -- (0,0) node[above right]{$(0,1)$} -- (-1,-0.6) node[below left]{$(1,k)$}
(0,0) -- (1,-0.6) node[below right] {$(1,q)$};
\filldraw (0,1.2) circle[radius=2pt] (0,0) circle[radius=2pt] (-1,-0.6) circle[radius=2pt] 
 (1,-0.6) circle[radius=2pt];
\end{scope}
\end{tikzpicture}

\end{tabular}
\end{center}

\caption{\footnotesize \label{fig surgeries}
Each column shows a different fixed point connected sum on the level of Coxeter polar data, 
the resulting gluing and its corresponding boundary tree. }
\end{figure}
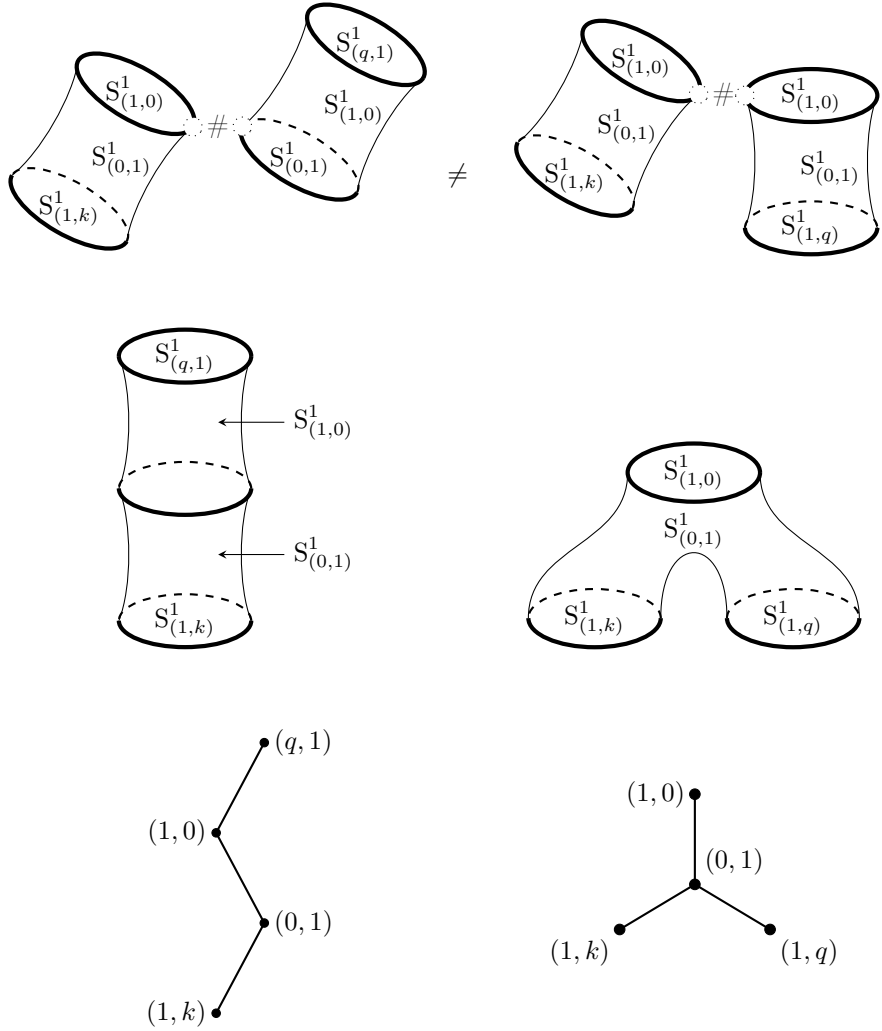

By means of fixed point connected sums  of the basic actions in \fref{fig basic boundary trees} 
we obtain all Coxeter polar $\T^2$ actions with orbit space $M^5/\T^2\cong\DD^3$.
\begin{prop} \label{prop T2 in M5 connected boundary}
 A Coxeter polar $\T^2$-manifold of dimension $5$ with orbit space homeomorphic to the $3$-disc and $N\geq 0$ 
 geodesics fixed by the action is either:
\begin{enumerate} 
 \item The factor-wise $\S^1\times \S^1$ polar action on $\mathbb{S}^4\times \mathbb{S}^1$, for $N=0$; 
 \item The linear $\T^2$ action on $\mathbb{S}^5$ of \eref{ex linear polar}, for $N=1$;
 \item The $\T^2$-actions on $\mathbb{S}^3\times_k \mathbb{S}^2$ described in \eref{ex T3 case main example} above, if $N=2$;
 \item Iterated $(N-2)$ fixed-point connected sums of $\mathbb{S}^3\times_k \mathbb{S}^2$, for different $k\in \ZZ$, in the case of $N\geq 3$.
\end{enumerate} \end{prop}

\begin{pf}
Notice that the number of fixed point connected components is the number of edges in the boundary tree,
and that we have already discussed the situation for $N\leq 2$. 
Observe as well that, if we consider two $\T^2$-manifolds $M_1$ and $M_2$ as above that are associated to connected boundary trees 
$\TT_i$, for $i=1,2$, and we perform an equivariant 
connected sum at respective fixed points, as in \pref{fixed point sum}, we will retrieve a new manifold with connected boundary tree given by 
the gluing along corresponding edges (the ones associated to the fixed geodesics containing those fixed points): 
$$ \TT = \TT_1 \bigcup_{e_1\sim e_2} \TT_2.$$
In particular, when we glue a linear sequence of three vertices along an edge we are adding one more vertex to the original tree. 
We can do this in any way we like, thus being able to ``grow'' all trees by this repeated operation, proving \mbox{part (4)}.
\qed
\end{pf}

\noindent \textit{Proof of \tref{THM T2 in M5}.} 
If we are given two Coxeter polar $\T^2$-manifolds of dimension five associated to polar data $(C_i,\T^2(C_i))$, $i=1,2$, 
we may perform a surgery along a regular orbit. By \dref{orbit sum},
the new chamber $C$ will be the connected sum through an interior point of the two given ones, $C= C_1 \# C_2$,  
and the new boundary will be the union of the preceding two, along with their unaffected marking. 
This can be translated as saying that the resulting boundary tree is the disjoint union of the two given ones.

In this way, we have shown how to retrieve all the admissible boundary trees or, equivalently, sets of Coxeter polar data, 
in terms of regular orbit sums of the actions listed in \pref{prop T2 in M5 connected boundary}.
\qed

\medskip
However, we would like to say more, and establish a reciprocal statement. 
\begin{thm} \label{T^2 in M^5 simply-connected conditions}
A Coxeter polar $\T^2$-manifold obtained from regular orbit sums of the actions in \pref{prop T2 in M5 connected boundary}
is simply-connected unless each summand is $\S^4\times \S^1$. 
In that case it will be simply-connected if and only if
the circle isotropy groups that occur generate the homology lattice of $\T^2$.
\end{thm}
\begin{pf}
Recall the decomposition of the manifold that defines a regular orbit sum
\begin{equation} \label{regular orbit sum decomp}
 (M_1 \#_{reg} M_2) \cong (M_1 - (\T^2\times \DD^3)) \cup_{\T^2\times \mathbb{S}^2} ( M_2 - (\T^2\times \DD^3)) 
\end{equation}
by using equation \eqref{orbit sum characterization} for the case of a regular torus orbit. 
We can apply Seifert-Van Kampen's formula on the previous decomposition to compute the fundamental group.
Notice first that removing a regular orbit does not affect the fundamental group of the manifold, because such orbits have \mbox{codimension $3$} in this case.
We then have:
\begin{equation} \label{VanKampen regular orbit sum}
 \pi_1
(M_1 \#_{reg} M_2) ~\cong~ [~\pi_1(M_1) \star \pi_1( M_2 ) ~]~/~ i_* \pi_1(\mathbb{S}^2\times \T^2).
\end{equation}
The manifolds we have to glue come from \pref{prop T2 in M5 connected boundary}, 
and they are all simply-connected with the exception of $\mathbb{S}^4\times \mathbb{S}^1$. 
For the latter, by the explicit description of the action we have that
\begin{equation} \label{SS^4 times SS^1}
(\mathbb{S}^4\times \mathbb{S}^1) - (\T^2 \times \DD^3) \cong (\mathbb{S}^4 - (\mathbb{S}^1 \times \DD^3)) \times \mathbb{S}^1 \cong \mathbb{S}^2 \times \DD^2 \times \mathbb{S}^1; 
\end{equation}
and the induced map coming from the inclusion of its boundary $\mathbb{S}^2\times \T^2$ is an epimorphism on the fundamental group level. 
Then, if we take regular orbit sums of $\mathbb{S}^4\times \mathbb{S}^1$ with an already simply-connected manifold 
we will obtain a new simply-connected manifold.
In particular, if the manifold admits fixed points, the boundary tree will have a connected component with at least one edge 
which corresponds to a simply-connected 
manifold we can begin with.

Otherwise, under the absence of fixed points, the boundary tree will have only isolated vertices marked by circle subgroups of slopes $(m_i,n_i)$ 
corresponding to regular orbit sums of blocks $(\mathbb{S}^4\times \mathbb{S}^1)_i$.
For each of these, we observe that the inclusion of a regular orbit induces an epimorphism in the fundamental groups, 
which leads to the following isomorphism
$$ \pi_1(~(\mathbb{S}^4\times \mathbb{S}^1)_i)= \pi_1(\T^2/\K_i) = ~\ZZ^2/<(m_i,n_i)> .$$

Fixing coordinates for $\T^2$ will simplify the computation of the fundamental groups along the gluing. 
For instance, if $M'$ is the regular orbit sum of two blocks $(\mathbb{S}^4\times \mathbb{S}^1)_i$, $i=1,2$, 
then \eqref{VanKampen regular orbit sum} implies that
\begin{equation*}
\begin{array}{ll}
\pi_1(M')	&\cong ~[~\ZZ^2/\!<\!(m_1,n_1)\!>] ~ \star~ [\ZZ^2/\!<\!(m_2,n_2)\!> ~] /~~\delta(\ZZ^2) \\
		&	\\
		&\cong ~~\ZZ^2 /<\!(m_1,n_1),(m_2,n_2)\!>,
\end{array}
\end{equation*}
where we have observed that the map corresponding to the inclusion of the boundary ${i_*}|_{\pi_1(\T^2)}$ 
is the diagonal map from $\ZZ^2$ onto each factor followed by the corresponding projections.
Moreover, the inclusion of a regular orbit generates the fundamental group of this new manifold $M'$.
In this way we may proceed inductively, where $\pi_1(M')\cong ~\ZZ^2/<\!(m_1,n_1),\cdots,(m_i,n_i)\!>$ is known,
and we want to compute the fundamental group of $M = (\mathbb{S}^4\times \mathbb{S}^1)_{i+1} ~\#_{reg} ~M'$. 
The inductive step is completely analogous, and so we have that the fundamental group of $M$ will be $\ZZ^2/<\!(m_1,n_1),\cdots,(m_{i+1},n_{i+1})\!>$.

In particular, the fundamental group will be trivial if and only if the isotropy subgroups generate the homology of $\T^2$, 
or cyclic if they generate a sub-lattice.
\qed
\end{pf}

\sect{Diffeomorphism Classification} \label{sect diffeo tori classif}
The goal of this section is to determine the diffeomorphism type of the underlying manifolds for the polar torus actions on 
simply-connected $5$-manifolds.
In the case of $\T^3$ this was done in \cite{Oh5} were it was proven that these are obtained as 
connected sums of $\mathbb{S}^5$, $\mathbb{S}^3\times \mathbb{S}^2$
and $\mathbb{S}^3\tilde{\times} \mathbb{S}^2$. This will also be true for $\S^1$ and $\T^2$ actions, as stated in \tref{THM diffeo tori in M5}.
We separate the two cases, addressing the simpler one of $\S^1$-actions first.
In both cases the proof relies on our explicit description of the actions, together with the
Barden-Smale classification of closed simply connected $5$-manifolds.
\begin{thm} \label{topological recognition cohom 4}
A compact simply-connected $5$-manifold admitting a polar circle action is diffeomorphic to a connected sum of
$\mathbb{S}^5$, $\mathbb{S}^3\times \mathbb{S}^2$, and $ \mathbb{S}^3\tilde{\times} \mathbb{S}^2$.
\end{thm}
\begin{pf}
Recall that these manifolds are reconstructed out of Coxeter polar data $(N,\S^1(N))$ 
given by a smooth simply-connected $4$-manifold $N$ with non-trivial boundary,
and a marking specifying the trivial subgroup along the interior and fixed points along the boundary.

The statement of the theorem follows from the classification of simply-connected closed $5$-manifolds (see \cite{Barden}), 
together with the claim that our manifolds have free second homology group $H_2(M;\ZZ)$.
To compute these groups we will use the Mayer-Vietoris sequence for a decomposition of $M$ given by 
\begin{equation} \label{MV space decomposition}
\partial N \times \S^1 \subset (\partial N \times \epsilon\DD^2) \cup (int(N) \times \S^1) = M,
\end{equation}
where  $\partial N \times \DD^2$ is a tubular neighborhood of $\partial N$ inside $M$, 
and $int(N) \times \S^1$ is the reconstructed manifold over the interior of $N$. 
Here, both the tubular neighborhood of $\partial N$ and the reconstructed manifold over $int(N)$
are identified due to the uniqueness properties of the reconstructed polar actions with prescribed data.
We then have the associated long exact sequence of homology with integer coefficients,
\begin{equation} \label{MV homology seq. circle M5}
H_{i}(\partial N \times \S^1) \xrightarrow{\alpha} H_{i}(\partial N) \oplus H_{i}(int(N) \times \S^1) 
\xrightarrow{\beta} H_{i}(M) 
\xrightarrow{\delta} H_{i-1}(\partial N \times \S^1)	.
\end{equation}
We have by hypothesis that $\pi_1(N)=0$, implying $H_1(N)=0$. 
Also notice that the interior of $N$ is a retract of it, and thus has the same homology.

By the K\"unneth formula, and thanks to the fact that $\S^1$ has free homology, 
we have:
$$\begin{array}{l}
H_1(\partial N \times \S^1) = H_1(\partial N) \oplus \ZZ^{\pi_0(\partial N)}, \\
      \\
H_2(\partial N \times \S^1) = H_2(\partial N) \oplus H_1(\partial N), \\
      \\
H_2(int( N) \times \S^1) = H_2(N) .
\end{array}$$
The map $\alpha$ in \eqref{MV homology seq. circle M5} gives isomorphisms between the corresponding terms in the $i$-th homology for $i=1,2$,
$$\alpha: H_i(\partial N) \subset H_i(\partial N \times \S^1) \xrightarrow{\sim} H_i(\partial N)  \cong H_i(\partial N \times \DD^2).$$
But there is another term coming from the common intersection in \mbox{$H_2(\partial N \times \S^1)$}, which is 
$H_1(\partial N) \otimes H_1(\S^1) \cong H_1(\partial N)$. In this case, its image via $\alpha$ will be trivial,
as the generator of $H_1(\S^1)$ is zero in $H_1(\DD^2)$, while the map from $H_1(\partial N)$ is trivial since $H_1(N)=0$.
Therefore, the kernel of $\delta$ in degree two is isomorphic to $H_2(int(N)\times \S^1) \cong H_2(N)$. 
It follows that $H_2(N)$ is free by the 
Poincar\'e duality for manifolds with boundary.

The image of $\delta$ is a submodule of $H_1(\partial N \times \S^1)$, and in particular of the summand 
$H_0(\partial N \times \S^1)$ which is free, as we have already seen how the term $H_1(\partial N)$ cancels out 
by applying itself isomorphically into its image via $\alpha$. We conclude that $H_2(M)$ is a free $\ZZ$-module.
\qed
\end{pf}

\begin{rem} \label{explicit out of data CIRCLES}
We may list the manifolds uniquely, with $\alpha=0,1$ and $n\in \NN_0$, as 
$$ \#_{\alpha} \mathbb{S}^3\tilde{\times} \mathbb{S}^2 \#_n \mathbb{S}^3\times \mathbb{S}^2,$$ 
where the empty connected sum refers to the $5$-sphere. 
These manifolds are uniquely identified by the rank of its second integer homology group, 
together with the property of being spin or not, which corresponds to $\alpha=0,1$ respectively.
This is explicitly determined by the Coxeter polar data. In fact, we can compute the rank of the second homology $H_2(M)$ as
$$H_2(M)= H_2(N)\oplus \ZZ^{|\pi_0(\partial N)| -1}.$$
Furthermore, $M$ will be spin if and only if $N$ is. 
First notice that, since $N$ retracts onto its interior, one being spin implies so for the other. 
Then, observe that, if $N$ is not spin, $M$ will not be spin as well,
for it contains the open subset $M_{reg}=int(N)\times \S^1$. 
On the other hand, if $N$ is spin we fix a spin structure which induces one on its interior and boundary.
We now can construct a spin structure on $M$ by gluing along the decomposition \eqref{MV space decomposition}.
Each piece has a product spin structure by choosing for $\S^1$ the unique spin structure that extends to $\DD^2$.
In this way, we ensure that the spin structures agree on the common 
open stripe where they are identified, allowing the gluing to be spin.
\end{rem}

\smallskip
For polar $\T^2$ actions, the identification of the underlying diffeomorphism types is analogous to the case of the circle actions, 
though significantly more involved. We rely on the general description of the actions achieved in \mbox{Section \ref{section Polar $T^2$ actions on $M^5$}.}
\begin{thm} \label{topological recognition T2 in M5}
 A compact simply-connected $5$-manifold admitting a polar effective $\T^2$-action is diffeomorphic to
 a connected sum of $\mathbb{S}^5,~\mathbb{S}^3\times \mathbb{S}^2,~\text{or}~\mathbb{S}^3\tilde{\times} \mathbb{S}^2$.
 \end{thm}
\begin{pf}
Analogously to the proof of \tref{topological recognition cohom 4} we will show that the second integer homology groups $H_2(M;\ZZ)$ are free.

Polar effective $\T^2$-actions on compact simply-connected $5$-manifolds were constructed by first forming connected sums of basic examples 
(this is the content of \pref{prop T2 in M5 connected boundary}), and later performing surgeries along normally embedded tori. 
For the computations that follow, it is useful to revert the order in which these operations are done. 
We will begin by taking regular orbits sums of the actions on $\mathbb{S}^4\times \mathbb{S}^1$ and $\mathbb{S}^5$ only, and later make fixed point sums with $\mathbb{S}^3\times_k \mathbb{S}^2$. 
This last step will only add free summands to the second homology, 
so we only need to compute the second homology for the initial regular orbit sums. 
From the decomposition \eqref{regular orbit sum decomp} 
we have the following Mayer-Vietoris long exact sequence for integer homology
 $$\cdots \to H_{2}(\mathbb{S}^2 \times \T^2 ) \xrightarrow{\alpha_2} H_{2}(M_1 - (\T^2\times \DD^3)) \oplus H_{2}(M_2 - (\T^2\times \DD^3)) \xrightarrow{\beta_2}  $$
\begin{equation} \label{mayer vietoris reg orbit sum homology}
\xrightarrow{\beta_2}  H_{2}(M_1 \#_{reg} M_2) \xrightarrow{\delta_2} H_{1}(\mathbb{S}^2 \times \T^2 ) \xrightarrow{\alpha_1} \cdots 	\hspace{3cm}
\end{equation}
In order to compute the homology of a regular orbit sum $H_{2}(M_1 \#_{reg} M_2)$, we first need to know what 
does each $M_i$ look like when a tubular neighborhood around a regular orbit has been removed. 
For $M_i=\mathbb{S}^4\times \mathbb{S}^1$ we have the equality \eqref{SS^4 times SS^1}, while for $\mathbb{S}^5$ we can make use of an auxiliary Mayer-Vietoris
computation on the decomposition given by a tube $\T^2\times \DD^3$ and its complement, and show that $ H_2(\mathbb{S}^5 - (\T^2\times \DD^3)) = \ZZ $. 
In both cases the second homology group is generated by the class of a $2$-sphere 
bounding a slice $\DD^3$ at a regular point. 
The group $H_{2}(M_1 \#_{reg} M_2)$ will then be computed as the kernel to $\alpha_1$ plus the image 
of $\beta_2$ from \eqref{mayer vietoris reg orbit sum homology}. 
This kernel is clearly free, as any subgroup of $\ZZ^2$, while the other summand requires
that we control the map $\alpha_2$ coming from the inclusion of the common intersection $\mathbb{S}^2\times \T^2$. 
Using the K\"unneth theorem for the product we have that 
$\ZZ^2 \cong H_2(\mathbb{S}^2\times \T^2)\cong \ZZ.[\mathbb{S}^2] \oplus \ZZ.[\T^2]$. 
As mentioned above, the class of the sphere in the common intersection $[\mathbb{S}^2] \in H_2(\T^2\times \mathbb{S}^2)$ 
generates, via $\alpha_2$, the second homology of the block that is being added, say $M_2$.
Then, what is left for the image of $\beta_2$ is just the other term corresponding to $M_1$, $H_2(M_1 -(\T^2\times \DD^3))$, 
when quotiented out by the image $\alpha_2(\ZZ.[\T^2])$.
This can be seen as the class of a regular orbit in the boundary of $(M_1 -(\T^2\times \DD^3))$. 

We claim that $\alpha_2([\T^2])\in H_2(M_1 -(\T^2\times \DD^3))$ is trivial since we can construct an embedded $3$-manifold with this torus for its boundary. 
For this, first consider the $3$-dimensional chamber associated to the $\T^2$-action in $M_1 -(\T^2\times \DD^3)$, 
which is just the chamber that we already got for $M_1$ minus an interior disc. 
Then, take a segment inside this chamber joining the point that corresponds to this regular orbit in the boundary together with 
a point in a face (with circle isotropy subgroup), arriving orthogonally to this last stratum. 
The orbit of this segment will be smoothly embedded in $M$, and diffeomorphic to $\DD^2\times \S^1$ with the required boundary, thus proving the claim.

In conclusion, if we can ensure that the piece \mbox{$M_1 -(\T^2\times \DD^3)$} has free second homology, this will also be true for the sum of it with a basic block,
\mbox{$M_1 \#_{reg} M_2$}.
In particular, it is true for the sum of two basic blocks. For the induction to continue, 
we must see that the complement of a tube around a regular orbit in the resulting manifold has free second homology.
A key observation in order to see this is that, instead of first gluing the pieces that define it
and then removing a tube along a regular orbit, 
we may directly glue the two pieces that define it along 
\mbox{$\T^2\times\DD^2$}, obtaining
$$ (M_1 \#_{reg} M_2) - (\T^2\times\DD^3) = (M_1 - (\T^2\times \DD^3)) \cup_{\T^2\times \DD^2} (M_2 - (\T^2\times \DD^3)).$$

\noindent Then, the homology can be computed via the corresponding Mayer-Vietoris sequence:
\begin{equation*} 
\begin{array}{l}
\cdots \to H_{2}(\T^2 \times \DD^2) \xrightarrow{\alpha_2} H_{2}(M_1 - (\T^2\times \DD^3)) \oplus H_{2}(M_2 - (\T^2\times \DD^3)) \xrightarrow{\beta_2} \\
\\
\xrightarrow{\beta_2}  H_{2}((M_1 \#_{reg} M_2) - (\T^2\times \DD^3)) \xrightarrow{\delta_2} H_{1}(\mathbb{S}^2 \times \T^2 ) \xrightarrow{\alpha_1} \cdots.
\end{array}
\end{equation*}
As before, $H_{2}(M_1 \#_{reg} M_2)$ can be decomposed into the clearly free kernel of $\alpha_1$ plus the image of $\beta_2$. 
The latter is isomorphic to the sum of \mbox{$H_{2}(M_i - (\T^2\times \DD^3))$}, for $i=1,2$, 
which are free by inductive hypothesis, once we
observe that the map $\alpha_2$ coming from the $[\T^2] \in H_2(\T^2\times \DD^2)$ 
is trivial by an analogous argument as in the previous claim.
\qed
\end{pf}

\begin{rem}
The diffeomorphism type can be determined explicitly from given Coxeter polar
data. Analogously to \rref{explicit out of data CIRCLES},
we may determine the rank of the second homology and characterize when it is spin.
In fact, we have a spin structure on any regular orbit sum of copies of $\mathbb{S}^4\times \mathbb{S}^1$ and $\mathbb{S}^5$. 
For this, it suffices to consider the given spin structure on $\mathbb{S}^5$, the unique spin structure on 
$\T^2\!\times \mathbb{S}^2\!\times \!(0,\epsilon)\subset \mathbb{S}^4\times \mathbb{S}^1$ that extends to $\mathbb{S}^4\times \DD^2$
and finally observe that the gluing map agrees with the identity along the tubular neighborhood of the regular orbit,
and is therefore spin.
On the other hand, 
if we perform a fixed point connected with an action on $\mathbb{S}^3\tilde{\times}\mathbb{S}^2$ the manifold will not be spin.
Hence, the resulting manifold will be spin if and only if it is constructed taking fixed point connected sums with $\mathbb{S}^3\times_k\mathbb{S}^2$ for even $k$ only.
\end{rem}

\sect{Applications to non-negative curvature \label{section Applications to non-negative curvature}}

The interest in positively and non-negatively curved compact manifolds
together with the abundance of symmetry for the relatively few known examples,
has led to a program of studying manifolds with positive or non-negative curvature on which a Lie group acts by isometries 
(see for example the surveys \cite{Grove,Wilking,Z}).

Remarkable results have been accomplished in low dimensions, 
where topological classifications are obtained for prescribed curvature and symmetry conditions.
For instance, compact simply-connected positively curved $4$-manifolds with continuous symmetry, i.e., 
invariant under an isometric circle action, 
are equivariantly diffeomorphic to $\mathbb{S}^4$ or $\CC\PP^2$ with a linear action (\cite{Hsiang-Kleiner,Fintushel,GroveWilking}). 
If we consider non-negatively curved ones, we have in addition $\mathbb{S}^2 \times \mathbb{S}^2$ and $\CC\PP^2 \# \pm \CC\PP^2$
(\cite{Kleiner_thesis,Searle-Yang}) and the circle action is induced by a biquotient action (\cite{GGK}).

In contrast, in dimension $5$ the topological question is open 
if we only consider the action of a circle, unless this action is fixed point homogeneous (see \cite{GGSpin}). 
Otherwise, we must strengthen the symmetry hypothesis as in \cite{Rong}, where it is proven that
a positively curved (closed simply-connected) $5$-manifold with an isometric \mbox{$\T^2$-action} is diffeomorphic to
the sphere $\mathbb{S}^5$.
The list grows if we allow the metric to be non-negatively curved to include the Wu manifold $\SU(3)/\SO(3)$, 
$\mathbb{S}^3\times \mathbb{S}^2$ and the non-trivial bundle $\mathbb{S}^3 \tilde{\times} \mathbb{S}^2$ (see \cite{GGS_nn_5d_almst_mxlsymrank}). 
However, equivariant classifications of non-negatively curved $5$-manifolds 
are only obtained in the case of an $\SU(2)$
action (see \cite{Simas}) or of a $\T^3$ action (see \cite{Oh5,GGK}).

Polar actions are expected to be rigid in this context. In fact, it was shown in \cite{FGT} that a 
polar action with cohomogeneity at least two and positively curved polar metric is equivariantly
diffeomorphic to a linear action on a CROSS (compact rank one symmetric space). 
However, the case of non-negative curvature is open in general.
Here, as an application of our classification results together with a closer examination of what the sections
are, we are able to classify polar actions with non-negative curvature in dimension $5$. 

\medskip
\noindent \textit{ Proof of \tref{THM non-negative curv}.} 
First observe that the actions listed in the theorem were given with an invariant polar metric of non-negative curvature. 
We can then focus on discarding the other cases that appear in our classification, 
using the fact that the sections have to inherit non-negative curvature since they are totally geodesic submanifolds. In particular, 
by a classic result of J. Cheeger and D. Gromoll, their fundamental groups 
have to admit a lattice $\ZZ^r$ of finite index. 

In general, the chamber $C$ of a polar manifold is a convex submanifold with totally geodesic boundary.
Moreover, if the curvature is non-negative, the proof of the Soul theorem \cite{Cheeger-Gromoll} implies that $C$ admits a soul $S$. 
Flowing along a gradient-like vector field adapted to the distance function to the soul, 
it follows that the chamber is homeomorphic to a disc bundle over $S$. 
The simply-connectedness of the manifold implies that of $C$ and hence of $S$. 

We separate the proof according to the Lie group being abelian or not. 
We begin with the $\S^1$ and $\T^2$ actions, and recall that the case of $\T^3$ actions was treated in \cite{Oh5,GGK}.

For a circle action, the soul can be a point, $\mathbb{S}^2$ or $\mathbb{S}^3$. In the first case we have $C= \DD^4$ which is recognized as a linear action on $\mathbb{S}^5$. 
For the last case, the chamber will be $\mathbb{S}^3\times I$, corresponding to the circle action on the second factor of $\mathbb{S}^3\times \mathbb{S}^2$. 
Finally, when the soul is $\mathbb{S}^2$ we have infinitely many possible disk bundles which gives rise
to actions on $\mathbb{S}^3\times_k \mathbb{S}^2$ as discussed in \eref{ex polar circle action non-neg}.

\smallskip
In the case of a $\T^2$ action, the soul of the chamber must be a point or $\mathbb{S}^2$, giving rise to $C=\DD^3$ or $\mathbb{S}^2\times I$, respectively.
We claim that only the basic actions on $\mathbb{S}^3 \times_k \mathbb{S}^2$, $\mathbb{S}^4\times \mathbb{S}^1$ and $\mathbb{S}^5$
of \fref{fig basic boundary trees}, and the linear $\T^2$ action on the first factor of $\mathbb{S}^3\times \mathbb{S}^2$ admit invariant non-negative curvature.
When the orbit space is $\DD^3$ the only other possible case is given by a repeated fixed point connected sum of $\mathbb{S}^3 \times_k \mathbb{S}^2$.
Notice that the actions on $\mathbb{S}^3 \times_k \mathbb{S}^2$ have $\Sigma = \mathbb{S}^2 \times \S^1$ or a non-orientable one, $\mathbb{S}^2\times \mathbb{S}^1/\tau_1$, 
as discussed in \eref{ex T3 case main example}.  
Then, the fundamental group of the section associated to a fixed point connected sum of at least two copies of $\mathbb{S}^3 \times_k \mathbb{S}^2$ will have
the free group on two generators as a subgroup and thus cannot admit non-negative curvature.
When the chamber is $\mathbb{S}^2\times I$ we have to perform one regular orbit sum between two of the manifolds listed in \pref{prop T2 in M5 connected boundary}. 
As follows from \pref{regular orbit sum section}, the fundamental group of such a gluing will be a free product 
of copies of the fundamental group of the summands
plus one additional generator for each 0-surgery. 
A case by case analysis shows that the only admissible fundamental group corresponds to both summands being $\mathbb{S}^4\times \mathbb{S}^1$.
Furthermore, by \tref{T^2 in M^5 simply-connected conditions} the gluing of two copies of $\mathbb{S}^4\times \mathbb{S}^1$ 
will add to a simply-connected manifold if and only if they are marked by complementary circle isotropy groups, 
i.e., circles generating $\T^2$ with trivial intersection. 
This is identified as the $\T^2$ action on $\mathbb{S}^3 \times \mathbb{S}^2$, coming from the action on the first factor only.

\smallskip
Among the polar actions by non-abelian Lie groups there are two families of actions 
listed in \tref{THM classification dim5} which are constructed as fixed point sums of 
linear polar actions on $\mathbb{S}^3\times \mathbb{S}^2$ by the groups $\SO(3)$ and $\SO(3)\times \S^1$. 
These basic actions admit $\mathbb{S}^1\times \mathbb{S}^2$ and $\S^1\times \S^1$ for their respective sections. 
In either case, a fixed point sum of two or more copies has a section which can not be non-negatively 
curved. 

Finally, among the $\SO(3)$ polar manifolds given by fixed point sums of the actions on the Brieskorn or Wu manifold $\#_{k}\B$ and $\# _{l}\W$, 
only the action on the Wu manifold alone supports non-negative curvature, because otherwise we have that $\Sigma=\#_n \T^2$, $n\geq 2$. 
\qed


\newpage

\bigskip
\bibliographystyle{alpha}

{\renewcommand{\baselinestretch}{1}
\hspace*{-20ex}\begin{tabbing}
IMPA - Rio de Janeiro - Brazil \\
E-mail address: fj.gozzi@gmail.com 
\end{tabbing}

}

\end{document}